\documentclass[regno, reqno,11pt]{amsart}
\usepackage{amsmath,amsthm,amssymb,amsfonts,epsfig,color,graphicx,enumerate,accents,subfigure,comment, esint,float,mathtools,appendix}
\usepackage[a4paper,margin=2.59truecm]{geometry}
\usepackage{graphicx}
\usepackage{enumitem}

\usepackage{hyperref}
\hypersetup{
  colorlinks   = true,    
  urlcolor     = blue,    
  linkcolor    = blue,    
  citecolor    = blue      
}

\usepackage{fancyhdr}
\DeclareMathOperator{\dive}{div}

\DeclareMathOperator{\diam}{diam}

\DeclareMathOperator{\dist}{dist}

\DeclareMathOperator{\Tr}{Tr}

\def\eps{{\varepsilon}}
\def\wto{\rightharpoonup}

\def\N{\mathbb{N}}

\def\R{\mathbb{R}}
\def\O{\Omega}
\def\o{\omega}

\def\vf{\varphi}

\def\GG{\mathcal G}

\def \< {\langle}
\def \>{\rangle}

\def\PP{\mathcal{P}}

\def\TT{\mathcal{T}}
\def\bf{{\Phi}}
\def\tt{\mathfrak{t}}

\newcommand{\be}{\begin{equation}}
\newcommand{\ee}{\end{equation}}

\newcommand{\bs}{\begin{split}}
\newcommand{\es}{\end{split}}

\newcommand{\red}{\textcolor{red}}

\newcommand{\LN}[2]{\| #1 \|_{L^N(#2)} }
\newcommand{\Linfty}[2]{\| #1 \|_{L^{\infty}(#2)} }
\renewcommand{\L}[3]{\| #2 \|_{L^{#1}(#3)}}
\renewcommand{\b}[1]{\bar{#1}}

\numberwithin{equation}{section}
\theoremstyle{plain}
\newtheorem{theorem}{Theorem}[section]
\newtheorem{lemma}[theorem]{Lemma}
\newtheorem{corollary}[theorem]{Corollary}
\newtheorem{proposition}[theorem]{Proposition}
\newtheorem{definition}[theorem]{Definition}
\theoremstyle{remark}
\newtheorem{remark}[theorem]{Remark}

\title[quasilinear transmission problem]{Existence and optimal regularity theory for weak solutions of free transmission problems of quasilinear type via Leray-Lions method.}

\author{Diego Moreira$^{*}$}
\address{$^*$Departamento de Matemática, Universidade Federal da Ceara(Fortaleza, Brazil)}
\email{$^*$dmoreira@mat.ufc.br}
\author{Harish Shrivastava$^\dagger$}
\address{$^\dagger$Tata Institute of Fundamental Researcher-Centre of of Applicable Mathematics}
\email{$^\dagger$harish21@tifrbng.res.in }

\begin{document}

\begin{abstract}
We study existence and regularity of weak solutions for the following PDE
$$
-\dive(A(x,u)|\nabla u|^{p-2}\nabla u) = f(x,u),\;\;\mbox{in $B_1$}.
$$
where $A(x,s) = A_+(x)\chi_{\{s>0\}}+A_-(x)\chi_{\{s\le 0\}}$ and $f(x,s) = f_+(x)\chi_{\{s>0\}}+f_-(x)\chi_{\{s\le 0\}}$. Under the ellipticity assumption that $\frac{1}{\mu}\le A_{\pm} \le \mu$, $A_{\pm}\in C(\O)$ and $f_{\pm}\in L^N(\O)$, we prove that under appropriate conditions the PDE above admits a weak solution in $W^{1,p}(B_1)$ which is also $C^{0,\alpha}_{loc}$ for every $\alpha\in (0,1)$ with precise estimates. Our methods relies on similar techniques as those developed by Caffarelli to treat viscosity solutions for fully non-linear PDEs (c.f. \cite{C89}) and also Leray-Lions method (c.f. \cite{BM92}, \cite{MT03}) to deal with compactness and existence of solutions. Moreover,  the $\TT_{a,b}$ operator (which was introduced in \cite{MS22}) also plays a key role in the regularity theory of solutions. 
\end{abstract}

\medskip 

\maketitle


\textbf{2010 Mathematics Subjects Classification:} 35R35, 35D30, 35J62 35B65, 35J60.

\tableofcontents
\section{Introduction}
Modelling the diffusion of various properties of materials, particularly the composite ones is a complicated process. In general, the diffusion depends on the values of the quantities representing what is being diffused, which in its turns depends on the material and physical properties of the medium. A typical model can be described as follows: let $\O$ be a given domain partitioned in mutually disjoint subsets $\{\O_i\}_{i=1}^k$ given by  
$$
\O_j = \O \setminus \big ( \cup_{i \neq j} \overline {\O_i }\big )
$$
and assume that inside $\O$ is prescribed the following PDE 
\be\label{basic form}
-\dive(F(x,u,\nabla u)) = f(x,u) 
\ee
with $F(x,s,\xi) = F_i (x,s,\xi)$ for $x\in \O_i$ and $f(x,s) = f_i(x,s)$ for $x\in \O_i$. In the formulation above, the diffusion processes may become discontinuous along the a-priori known boundaries $\partial \O_i$. This class of problems is called transmission problems and have been widely studied in \cite{C57}, \cite{C59}, \cite{C59.2}, \cite{LS55}, \cite{OO61}, \cite{MS60}, \cite{S56} etc. We refer the reader to introduction chapter of \cite{M10} for a detailed historical account on developments in transmission problems.

In this paper, we have chosen another model of the form \eqref{basic form}, namely, where $F$ and $f$ have a jump discontinuity across the level set $\{u=0\}$.  This means that $F$ and $f$ are given by
$$
\mbox{$F(x,s,\xi) = F_{\pm} (x,s,\xi) \chi_{\{s^{\pm} >0\}}$, $f(x,s) = f_\pm(x,s)\chi_{\{s^{\pm} >0\}}$.}
$$
The discontinuity represents the transition from one composite material to another. The class of problems mentioned above are generally called \textit{free} transmission problems, since the boundaries along which discontinuity takes place is a-priori unknown.

We study the weak solutions of the following PDE ($1<p<\infty$),
\be\label{PDE}
-\dive(A(x,u)|\nabla u|^{p-2}\nabla u) = f(x,u)\mbox{ in $\O\subset \R^N$}.
\ee
where $A(x,s) = A_+(x)\chi_{\{s>0\}}+A_-(x)\chi_{\{s\le 0\}}$ and $f(x,s) = f_+(x)\chi_{\{s>0\}}+f_-(x)\chi_{\{s\le 0\}}$. The following assumptions are enforced throughout the paper
\begin{enumerate}[label = \textbf{H\arabic*}.]
  \setcounter{enumi}{-1}
  \item \label{H0} $\O \subset \R^N$, is open and bounded.
\item\label{H1} (\textbf{Continuity}) $A_{\pm}$ are  continuous in $\O$.
\item\label{H2} $f\in L^N(\O)$.
\item\label{H3} (\textbf{Ellipticity}) There exists $0<\mu<1$ such that for every $x\in \O$ we have $\mu  \le A_{\pm}(x) \le  \frac{1}{\mu}$.
\end{enumerate}

Recent developments in the case of free transmission problems can be found in \cite{TA15}, \cite{MS22}, \cite{DH21}, \cite{CSS21}, \cite{PA22} and \cite{MS22.2}. In \cite{MS22}, we have shown optimal regularity for solutions to a variational free transmission problem (for the case $p=2$). Both works \cite{TA15} and \cite{MS22} assume that the coefficients $A_{\pm}$ are only continuous. In \cite{MS22.2}, we assumed Hölder regularity of coefficients and showed that under appropriate boundary data, free boundary and fixed boundary touch each other tangentially. We point out that in \cite{KLS19} the equation \eqref{PDE} is treated in the case $p=2$ with $a_{\pm}(x)$ Hölder continuous and zero right hand side. In this case, local Lipschitz regularity is obtained. This paper deal with more general situation. Here, we consider the general case $p\in (1,\infty)$ with only continuous coefficients $a_{\pm}$ and a non-zero right hand side in $L^N(\O)$. Under these conditions, it is already known that the solutions are no longer locally Lipschitz (c.f. \cite{JMS09}) even in the case $a_+ = a_-$ (i.e. no transmission case) and $f_{\pm} =0$. So, our main result here states that, under the assumptions mentioned before solutions to \eqref{PDE} are $C^{0,1^-}$ with (precise) estimates. In \cite{CKS21}, variational formulation with different exponent in each phase is discussed. We refer to the references therein for more recent developments on the subject. 

This paper can be roughly divided into two macro parts. In the first one, we prove the existence of weak solutions to the PDE \eqref{PDE} by mollifying the problem via a parameter $\eps>0$ and then passing to the limit $\eps\to 0$. {We develop a strategy similar to the Leray-Lions method that accounts for almost everywhere convergence of the gradients. In fact, the behaviour of weak convergence under non-linearities is a delicate issue, see for instance \cite{MT03}. Here, it is indeed the place where the quasilinear version of our problem, (i.e. $p\neq 2$) brings substantial new difficulties (compare with variational ``linear" version (i.e. $p=2$) in our recent work \cite{MS22}).

Similar issue is also addressed in the by now classical paper by Boccardo and Murat \cite{BM92}.} In Section \ref{compactness}, we obtain a compactness lemma (c.f. Proposition \ref{step 1}, Proposition \ref{step 1.5}) which essentially says that small perturbations in $A_{\pm}$ and $f_{\pm}$ imply that weak solutions to \eqref{PDE} are as close as we wish to regular profiles. {The Leray-Lions method is once more used in the proof of compactness lemma. This is somehow expected, once our method to prove the existence of weak solutions  (c.f. Proposition \ref{existence}) is based on sequences of approximating  problems parametrised by $\eps>0$.}  

For the optimal growth rate of weak solutions along the zero level set, we implement an approximation theory similar to the one developed by L. Caffarelli in the seminal paper \cite{Ca89} to treat regularity theory of viscosity  solutions to fully non-linear PDEs. Here, there are mainly two main steps. In the first one, we prove the proximity of weak solutions to regular profiles by compactness argument in small scales (c.f. Proposition \ref{step 1} and Proposition \ref{step 1.5}). In our case, the regular profiles inherit regularity from the weak solutions of \eqref{PDE} with $A_{\pm}$ being constants, via $\TT_{a,b}$ operator discussed in the sequel. The second one is to reduce the problem to the so called ``small regime configuration" (c.f. Lemma \ref{step 2}) via the scaling invariance of the problem under the appropriate regularity assumptions of the data.

Just like \cite{MS22}, the primary tool utilized in proving optimal regularity of weak solutions for \eqref{PDE} is the $\TT_{a,b}$ operator. A detailed discussion on $\TT_{a,b}$ operator can be found in \cite[Section 3]{MS22}. We point out that our methods in the proofs of existence and regularity theory are different from the ones in existing literature. In particular, the techniques used in this paper allows to obtain precise estimates in any dimension $N\ge 2$.

\begin{definition}
$u\in W^{1,p}(\O)$ is called a weak solution of \eqref{PDE} if for all $\vf \in W_0^{1,p}(\O)$ the following equation is true
\be\label{weak solution}
\int_{\O} \Big ( A(x,u)|\nabla u|^{p-2}\nabla u\cdot \nabla \vf \Big )\,dx = \int_{\O} f(x,u) \vf\,dx
\ee
\end{definition}
\subsection{Main results}
Below, we list the main results of this paper. We highlight here that no conditions relating $p$ and $N$ are imposed. In particular, the results obtained in the paper hold  in any dimension $N\ge 2$ and $p\in (1,\infty)$ as far as regularity is concerned (some exceptions appear in the existence part).  As mentioned before, our results also encompass sharp regularity theory for classical (i.e. no transmission between the phases) quasilinear elliptic PDE in divergence form with continuous coefficients, bringing precise estimates. The proof of the theorem below can be found in Appendix \ref{Appendix A}. For related results one can also check \cite{EVT11}. 
\\
\\
In the sequel, $p\in (1,\infty)$ and $N\ge 2$ are conditions enforced throughout the entire paper. Exceptions (only in the existence part w.r.t. $p$) will be mentioned explicitly.

\begin{theorem}\label{supporting lemma}
Assume $u$ is a bounded weak solution of  
\be\label{simple PDE}
-\dive(A(x)|\nabla u|^{p-2}\nabla u) = f \;\;\mbox{in $B_1$}
\ee
with $A\in C(B_1)$, $\mu\le A \le \frac{1}{\mu}$ for some $\mu>0$, and $f\in L^N(B_1)$. Then, for any $0<\alpha<1$, $u\in C^{0,\alpha}_{loc}(B_1)$. In particular, for every $r<1$ we have the following estimate.
\be\label{estimates}
\|u\|_{C^{0,\alpha}(B_r)} \le \frac{C(N, p, \mu , \alpha,\o_{A,B_{r^*}})}{(1-r)^\alpha} \big ( \Linfty{u}{B_{1}} + \LN{f}{B_{1}}^{\frac{1}{p-1}} \big )
\ee
where $r^*:=\frac{1+r}{2}$ and $\o_{A,B_{r^*}}$ modulus of continuity of $A$ in $B_{r^*}$.
\end{theorem}
\begin{remark}
We would like to thank Prof. Giuseppe Mingione for bringing the paper \cite{KM12} to our attention. Theorem \ref{supporting lemma} has also been obtained by different methods with slightly different presentation.
\end{remark}
Main results of the paper are the following
\begin{proposition}[Existence of weak solutions]\label{existence}
If $f_{\pm}\in L^{\infty}(\O)$ and $p\ge 2$, there exists a weak solution of the PDE \eqref{PDE} in $W^{1,p}(\O)$ and all the weak solutions are locally bounded.
\end{proposition}

\begin{remark}
{The assumptions $f_{\pm}\in L^{\infty}(\O)$ do not appear to be optimal for existence of weak solutions of PDE \eqref{PDE}. Our methods in the proof of Proposition \ref{existence} seems to be working only for the case  $p\ge 2$.  We leave as an interesting open question, the existence of weak solutions to \eqref{PDE} for more general class of RHS and $1<p<2$. We point out however, that in the upcoming sections, we deal with the regularity theory for weak solutions under the more general assumptions that $f_{\pm}\in L^N(\O)$ and $1<p<\infty$. }
\end{remark}
{We already know that weak solutions of \eqref{PDE} are locally Hölder continuous in $\O$ (c.f. \cite[Chapter 10, Theorem 3.1]{ED09}). That is, any weak solution of \eqref{PDE} $u\in C^{0,\beta_0}_{loc}(\O)$ for some $\beta_0(p, \mu, f_{\pm}, N)$. In this paper we prove the local Hölder continuity for a weak solution of \eqref{PDE}, for any exponent $\alpha\in (0,1)$ along with estimates. }
\begin{theorem}\label{main result}
Assume $u$ is a bounded weak solution of \eqref{PDE} in $B_1$ satisfying the assumptions \eqref{H0}-\eqref{H3}. Then $u\in C^{0,1^-}_{loc}(B_1)$ with estimates. More precisely, for any $0<\alpha<1$, $u\in C^{0,\alpha}_{loc}(B_1)$ and for every $r<1$ we have the following estimates.
\be\label{estimates}
\|u\|_{C^{0,\alpha}(B_r)} \le \frac{C(N, p, \mu ,\alpha, \o_{A_{\pm},r^*})}{(1-r)^\alpha} \Bigg ( \Linfty{u}{B_{1}} + \LN{f_+}{B_{1}}^{\frac{1}{p-1}} + \LN{f_-}{B_{1}}^{\frac{1}{p-1}} \Bigg )
\ee
where $r^*:= \frac{1+r}{2}$  and $\o_{A_{\pm},{r^*}}$ is maximum of modulus of continuity of $A_+$ and $A_-$ in $\overline{B}_{r^*}$.
\end{theorem}

\begin{remark}
In fact, there is a slightly more general estimate that holds in the Theorem \ref{main result} above. Namely for any $D\subset \subset B_1$, we can prove the following estimate
\be\label{general est...}
[u]_{C^{\alpha}(D)} \le \frac{C(N,p,\mu,\alpha,\o_{A_{\pm},{d}})}{\dist(D,\partial B_1)^{\alpha}}\left ( \Linfty{u}{B_{1}}+ \LN{f_+}{B_{1}}^{\frac{1}{p-1}} + \LN{f_-}{B_{1}}^{\frac{1}{p-1}} \right ).
\ee
where $D\subset \subset B_{d} \subset \subset B_1$ and $d=1-\frac{\dist(D,\partial B_1)}{2}$ and $\o_{A_{\pm},{d}}$ is maximum of modulus of continuity of $A_+$ and $A_-$ in $\overline{B}_{d}$. A similar remark is true regarding Theorem \ref{supporting lemma} with the appropriate straightforward changes.
 
\end{remark}

\section{Preliminary definitions and supporting lemmas}
We set the following notation for $\vf \in W^{1,p}(D)$
\be\label{affine sobolev}
W_{\vf}^{1,p}(D) : = \Big \{ v\in W^{1,p}(D)\;:\: v-\vf \in W_0^{1,p}(D) \Big \}.
\ee

\begin{definition}$($\cite[Definition 3.1]{MS22}$)$\label{Tab}
Let $a,b>0$, $p\in [1,\infty)$ and $D$ be an open Lipschitz set, we define $\TT_{a,b}:W^{1,p}(D)\to W^{1,p}(D)$ as follows
$$
\TT_{a,b}(v) := a v^+\,-\,bv^-.
$$
We define the $\TT_{a,b}$-operator also on the boundary level acting on $L^p(\partial D)$ to be \\  $\TT^{\partial D}_{a,b}:L^p(\partial D) \to L^p(\partial D)$ given by 
$$
\TT^{\partial D}_{a,b}(\varphi) := a \varphi^+\,-\,b\varphi^-.
$$
\end{definition}

One can easily verify the following (c.f. \cite[Section 3]{MS22}) for all $q\in \R$
\be\label{gradient relation}
\begin{split}
|\nabla (\TT_{a,b}(u))|^{q}  &= |\nabla (au^+)|^{q}+ |\nabla (bu^-)|^{q}.\\
|\nabla (\TT_{a,b}(u))|^{q-2}\nabla (\TT_{a,b}(u)) &= |\nabla (au^+)|^{q-2}\nabla (au^+) - |\nabla (bu^-)|^{q-2}\nabla (bu^-)\\
&=a^{q-1}|\nabla u^+|^{q-2}\nabla u^+ - b^{q-1}|\nabla u^-|^{q-2}\nabla u^-.
\end{split}
\ee
As we have also seen in \cite{MS22}, the $\TT_{a,b}$ operator preserves any regularity up to Lipschitz. This means
\begin{lemma}\label{preservation of regularity}
For a function $u:D\to \R$ and $\alpha \in [0,1]$, the following holds true
$$
u\in C^{0,\alpha}(D) \iff \TT_{a,b}(u) \in C^{0,\alpha}(D).
$$
Moreover, the following estimate holds 
$$
[u]_{C^{0,\alpha}(D)} \le \frac{1}{\min(a,b)}[\TT_{a,b}(u)]_{C^{0,\alpha}(D)}.
$$
\end{lemma}
\begin{proof}
We prove one way of the implication, the other way follows by similar arguments. Let $\TT_{a,b}(u) \in C^{0,\alpha}(D)$, then for $x,y \in D$ , the following holds
\[
\begin{split}
x,y\in \{u>0\} &\implies |u(x)-u(y)| = \frac{1}{a}|au(x) - au(y)| =\frac{1}{a}|\TT_{a,b}(u)(x) -\TT_{a,b}(u)(y)|\le \frac{[\TT_{a,b}(u)]_{C^{0,\alpha}(D)}}{a}|x-y|^\alpha\\
x,y\in \{u\le0\}& \implies |u(x)-u(y)| = \frac{1}{b}|bu(x) - bu(y)| =\frac{1}{b}|\TT_{a,b}(u)(y) -\TT_{a,b}(u)(x)|\le \frac{[\TT_{a,b}(u)]_{C^{0,\alpha}(D)}}{b}|x-y|^\alpha\\
\end{split}
\]

\[
\begin{split}
x\in \{u>0\},y\in \{u\le 0\} \implies |u(x)-u(y)| =\Big  |\frac{1}{a}au(x) -\frac{1}{b} b u(y) \Big | &=\Big |\frac{1}{a}\TT_{a,b}(u)(x) -\frac{1}{b} \TT_{a,b}(u)(y) \Big |\\
&\le \frac{1}{\min(a,b)} [\TT_{a,b}(u)]_{C^{0,\alpha}(D)} |x-y|^\alpha\\
\end{split}
\]
The above three implications imply that 
$$
[u]_{C^{0,\alpha}(D)} \le \frac{1}{\min(a,b)}[\TT_{a,b}(u)]_{C^{0,\alpha}(D)}.
$$
\end{proof}
We present another important supporting lemma which comes handy while invoking the Leray-Lions method in upcoming theorems. Before that we quote a well known result (c.f \cite[Proposition 17.3]{CM09}): For $p\in (1, \infty)$ and $\xi, \zeta \in \R^N$ the following holds 
\be\label{convex identity}
\begin{split}
  (|\xi|^{p-2} \xi - |\zeta|^{p-2} \zeta )\cdot (\xi-\zeta) &\ge 0 \mbox{ more precisely,}\\
 (|\xi|^{p-2} \xi - |\zeta|^{p-2} \zeta )\cdot (\xi-\zeta)&\ge \begin{cases}
 |\xi-\zeta|^2 (|\xi|+|\zeta|)^{p-2} \ge |\xi-\zeta|^p \qquad \qquad  \mbox{        if $p\ge 2,$}\\
 |\xi-\zeta|^2 (|\xi|+|\zeta|)^{p-2}\;\;\;\;\;\;\;\;\;\;\;\;\;\;\;\;\;\;\;\;\;\;\;\;\;\;\;\;\;\; \mbox{              if $1<p< 2$}.
 \end{cases}
 \end{split}
\ee
  
\begin{lemma}\label{claim}
Let $v_k,v \in L^p(\O)^N$ be two sequences and $\GG_k$ be defined as 
$$
\GG_k := (|v_k|^{p-2}v_k - |v|^{p-2}v)\cdot (v_k-v).
$$ 
The following implication holds true
$$
\GG_k \to 0 \mbox{ a.e. in $\O$} \implies \;\;(v_k -v)\to 0 \mbox{ a.e. in $\O$}.
$$
\end{lemma}
\begin{proof}
To prove Lemma \ref{claim}, let us denote the sets $Z$ and $Z_0$ as follows
\[
\begin{split}
Z &: = \Big \{ x\in \O: \lim_{k\to \infty} \GG_k(x) \neq 0, \mbox{ or $\lim_{k\to \infty} \GG_k(x)$ does not exist}  \Big \}\\
Z_0 &:= \Big \{ x\in \O: |v| = \infty \Big \}\\
C_0&:=\O\setminus (Z\cup Z_0).
\end{split}
\]
Since $v\in L^p(\O)^N$, $|v|$ is finite almost everywhere in its domain and since $\GG_k\to 0$ a.e. in $\O$, therefore $|Z\cup Z_0|=0$. For the case $p\ge 2$, from \eqref{convex identity} we have 
$$
0\le |v_k-v|^p \le ( |v_k|^{p-2}v_k  - |v|^{p-2} v\big ) \cdot  (v_k-v)=\GG_k  \to 0 \mbox{ pointwise in $C_0$}.
$$
Hence, we obtain 
\be\label{c0}
p\ge 2 \implies |v_k-v|^p \to 0 \implies v_k \to v \mbox{ a.e. in $\O$}.
\ee
Thus the claim in Lemma \ref{claim} is true when $p\ge2$. In order to prove the claim for the case $1<p<2$. Let $x_0\in C_0:=\O\setminus (Z\cup Z_0)$, we consider three cases:
\begin{enumerate}[label = \textbf{(\alph*)}]
\item \label{a} $1<p<2$, $|v(x_0)|=0$.
\item \label{b} $1<p<2$, $|v(x_0)| >0$ and the sequence $v_k(x_0)$ is bounded.
\item \label{c} $1<p<2$, $|v(x_0)| >0$ and the sequence $v_k(x_0)$ is unbounded.
\end{enumerate}
From \eqref{convex identity} we have 
\be\label{p<2}
\begin{split}
 0\le \frac{|v_k -v|^2}{(|v_k| + |v|)^{2-p}} &\le ( |v_k|^{p-2}v_k  - |v|^{p-2} v\big ) \cdot  (v_k-v)=\GG_k \to 0 \mbox{ pointwise in $C_0$}\\
 &\implies \frac{|v_k -v|^2}{(|v_k| + |v|)^{2-p}} \to 0 \mbox{ pointwise in $C_0$.}
 \end{split}
\ee

In the case when $|v(x_0)|=0$, from \eqref{p<2} we have
\be\label{c1}
\mbox{Case \ref{a}} \implies |v_k(x_0)|^p \to 0 \implies |v_k(x_0)| \to  |v(x_0)|=0.
\ee
For the case \ref{b}, let $M_0 = \limsup_{k\to \infty}|v_k(x_0)|$, from \eqref{p<2}, for $k$ sufficiently large we have  
\be\label{c2}
\begin{split}
\mbox{Case \ref{b}}&\implies  0\le \frac{|v_k(x_0) -v(x_0)|^2}{(M_0 + |v(x_0)|)^{2-p}} \le \frac{|v_k(x_0) -v(x_0)|^2}{(|v_k(x_0)| + |v(x_0)|)^{2-p}} \to 0\\
& \implies v_k(x_0) \to v(x_0)\mbox { in $C_0$}.
\end{split}
\ee
If $v_k(x_0)$ is an unbounded sequence, then for $k$ sufficiently large we have 
\[
\begin{split}
\frac{1}{2}|v_k(x_0)|   \le |v_k(x_0) -v(x_0)| \mbox{ and }
(|v_k(x_0)| + |v(x_0)|) \le 2|v_k(x_0)|.
\end{split}
\]
From \eqref{p<2}, the equation above implies 
\be\label{c3}
\mbox{Case \ref{c}} \implies 0 \le |v_k(x_0)|^p \le \frac{|v_k(x_0) -v(x_0)|^2}{(|v_k(x_0)| + |v(x_0)|)^{2-p}} \to 0.
\ee
The implication above in \eqref{c3} is clearly a contradiction. Therefore, the case \ref{c} can not hold. From \eqref{c1} and \eqref{c2}, we prove the claim when $1<p<2$. This finishes the proof of Lemma \ref{claim}.
\end{proof}

\section{Existence of weak solution (Proof of Proposition \ref{existence})}
{
\begin{proof}[Proof of Proposition \ref{existence}]
Our strategy to prove the existence of a weak solution for \eqref{PDE} is to mollify the PDE \eqref{PDE} by a parameter $\eps>0$ to \eqref{eps PDE} with corresponding weak solution $u_{\eps}$. Then we pass to the limit $\eps\to 0$ and show that the limit(s) $\lim_{\eps\to 0} u_{\eps}$ is (are) weak solution(s) to the desired PDE \eqref{PDE}. We proceed step by step.
\\
\\
\textbf{STEP 1. Construction of a regularised problem.} 
\\
\\
In order to mollify the broken coefficients, let $\eps>0$ and we define the function $\psi^+$ and accordingly $\psi^+_\eps$.
$$
\psi^+(t) := \begin{cases}
1\;\;\mbox{for $t\ge1$}\\
t\;\;\mbox{for $0<t<1$}\\
0\;\;\mbox{for $t\le 0$}.
\end{cases}, \;\;\psi^+_{\eps}(t) := \psi^+(\frac{t}{\eps}).
$$
Accordingly, $\Psi^+_{\eps}$ is defined as
$$
\Psi^+_{\eps}(t):= \int_{-\infty}^t \psi^+_{\eps}(s)\,ds.
$$

also we define the functions $\psi^-_{\eps}$ and $\Psi^-_{\eps}$ as follows:
$$
\Psi^-_{\eps}(t) := \Psi^+_{\eps}(t) - t\mbox{ and } \psi^-_{\eps}(t):= \Big | \frac{d}{dt}\Psi^-_\eps(t) \Big | = 1-\psi^+_{\eps}(t).
$$
 \begin{figure}[H]\label{fig1}
 \floatplacement{figure}{H}
\includegraphics[width=0.7\textwidth]{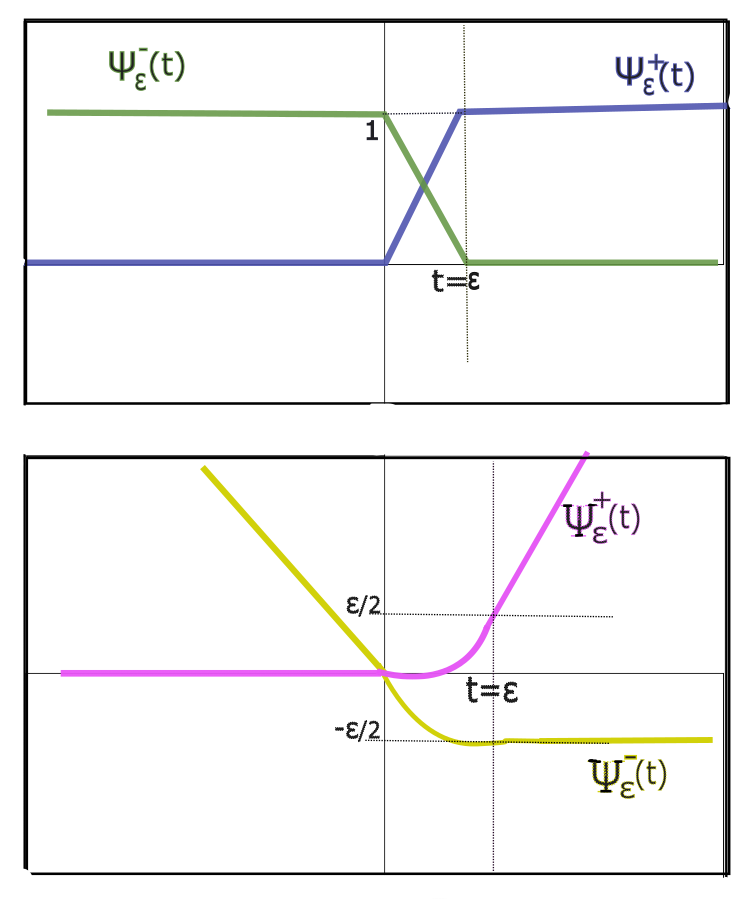}
\caption{(graphs of mollifiers)}
\end{figure}
Sometimes we will use alternate notations for chatacteristic functions for $\{s>0\}$, $\{s\le 0\}$
\[
\begin{split}
\psi_0^+(s)&:= \chi_{\{s>0\}}\\
\psi_0^-(s)&:= \chi_{\{s\le 0\}}
\end{split}
\]
We observe that 
\be\label{convg to char}
\psi_\eps^\pm \to \psi_0^\pm \mbox{ a.e. in $\R$.}
\ee
We define the mollified coefficient as follows
\be\label{def1}
\begin{split}
A_{\eps}(x,s):= A_+(x) ({\psi^+_{\eps}} (s) )^{p-1}+A_-(x)({\psi^-_{\eps}} (s) )^{p-1}\\
f_{\eps}(x,s):= f_+(x) \psi^+_{\eps} (s)+ f_-(x)\psi^-_{\eps}(s).
\end{split}
\ee
The mapping $(x,s)\to A_{\eps}(x,s)$ is continuous in $x$ and Lipschitz in $s$, therefore, for a given $g\in W^{1,p}(\O)$, there exists a weak solution $u_\eps$ for the following PDE
\be\label{eps PDE}\tag{$\PP_{\eps}$}
\begin{cases}
\dive(A_{\eps}(x,u_{\eps}) |\nabla u_{\eps}|^{p-2}\nabla u_{\eps})  = f_{\eps}(x,u_{\eps}) \;\;\mbox{in $\O$},\\
u_{\eps}-g \in W_0^{1,p}(\O).
\end{cases}
\ee
That is, $u_\eps\in W_g^{1,p}(\O)$ be such that for all $\vf \in W_0^{1,p}(\O)$ 
\be\label{eps PDE weak}
\begin{split}
\int_{\O}\Big(  A_{\eps} (x,u_{\eps}) |\nabla u_{\eps}|^{p-2} \nabla u_{\eps} \cdot \nabla \vf \Big )\,dx = 
 \int_{\O} f_{\eps}(x,u_{\eps}) \vf \,dx 
\end{split}
\ee
\begin{remark}
{We can verify that the PDE \eqref{eps PDE} satisfies the assumptions in (3.5), (3.6), (6.1), (6.2) and (6.3) in \cite{MZ97} when $f_{\pm}\in L^{\infty}(\O)$. The existence of a weak solution $u_\eps$ (when $f_{\pm}\in L^{\infty}(\O)$) follows from \cite[Theorem 6.12]{MZ97}. However, in upcoming sections, we have proved the regularity theory for weak solutions under the assumption $f_{\pm}\in L^N(\O)$. Authors leave the question regarding the existence of weak solutions to \eqref{eps PDE} for more general class of RHS as an open question.}
\end{remark}
\textbf{Step 2: Uniform $W^{1,p}(\O)$ bounds on $\{ u_{\eps}\}_{\eps>0}$}.
\\
\\
Since $(u_\eps-g) \in W_0^{1,p}(\O)$, by the definition of weak solution of \eqref{eps PDE}, we have the following
\[
\begin{split}
\int_{\O}\Big(  A_{\eps} (x,u_{\eps}) |\nabla u_{\eps}|^{p-2} \nabla u_{\eps} \cdot \nabla (u_{\eps}-g)\Big )\,dx  &= \int_{\O} f_{\eps}(x,u_{\eps}) (u_{\eps}-g)\,dx.\\
\end{split}
\]
Splitting and rearranging the terms we obtain
\be\label{eqn1.1}
\begin{split}
\int_{\O}  A_{\eps} (x,u_{\eps}) |\nabla u_{\eps}|^{p}\,dx = \int_{\O} f_{\eps}(x,u_{\eps}) (u_{\eps}-g)\,dx + \int_{\O}  A_{\eps} (x,u_{\eps}) |\nabla u_{\eps}|^{p-2}\nabla u_{\eps}\cdot \nabla g\,dx
\end{split}
\ee
We note that, from the definition of $A_{\eps}(x,s)$ and ellipticity of $A_{\pm}$ (c.f. \ref{H3}), we have 
\be\label{Aeps bounded}
2\mu \le A_{\eps}(x,s)\le \frac{2}{\mu}.
\ee 
We observe that if $p<N$, by Gagliardo-Nirenberg-Sobolev inequality, we have (here $F:= |f_+|+|f_-|$)
\be\label{controlling RHS 1}
\begin{split}
\int_{\O} \big | f_{\eps}(x,u_{\eps}) (u_{\eps}-g) \big | \,dx  &\le  \L{{p^*}'}{F}{\O}  \L{p^*}{u_{\eps}- g}{\O}\\
&\le C(\O) \L{{p^*}'}{F}{\O}  \L{p}{\nabla (u_{\eps}-g)}{\O}\\
&\le C(N,p,\O) \L{N}{F}{\O}  \L{p}{\nabla (u_{\eps}-g)}{\O} 
\end{split}
\ee
If $p\ge N$, then $u_{\eps}-g \in L^q(\O)$ for every $1<q<\infty$, therefore
\be\label{controlling RHS 2}
\begin{split}
\int_{\O} \big | f_{\eps}(x,u_{\eps}) (u_{\eps}-g) \big | \,dx  &\le  \L{N}{F}{\O}  \L{N'}{u_{\eps}- g}{\O}\\
&\le C(N,p,\O) \L{N}{F}{\O}  \|{u_{\eps}- g}\|_{W^{1,p}(\O)}\\
&\le C(N,p,\O)\L{N}{F}{\O}  \|\nabla (u_{\eps}- g)\|_{L^{p}(\O)}
\end{split}
\ee
Since $|\nabla u_\eps|\in L^p(\O) \implies |\nabla u_\eps|^{p-1} \in L^{p'}(\O)$, from Cauchy-Schwarz inequality we have 
\be \label{controlling LHS}
\begin{split}
\int_{\O}  A_{\eps} (x,u_{\eps}) |\nabla u_{\eps}|^{p-2}\nabla u_{\eps}\cdot \nabla g\,dx & \le C(\mu) \int_{\O} |\nabla u_{\eps}|^{p-1}|\nabla g|\,dx\\
&\le C(\mu) \L{p}{\nabla u_{\eps}}{\O}^{p-1}\cdot\L{p}{g}{\O}.
\end{split}
\ee
Plugging in the estimates \eqref{controlling RHS 1}, \eqref{controlling RHS 2} and \eqref{controlling LHS} in  \eqref{Aeps bounded} on rest of the terms in \eqref{eqn1.1} we obtain 
\[
\begin{split}
\int_{\O}|\nabla u_{\eps}|^p\,dx  &\le  C(\mu) \Bigg [ C(N,p,\O) \L{N}{F}{\O}  \L{p}{\nabla (u_{\eps}-g)}{\O}  + \L{p}{\nabla u_{\eps}}{\O}^{p-1}\cdot\L{p}{g}{\O} \Bigg ]\\
&\le C(\mu) \Bigg [ C(N,p,\O) \L{N}{F}{\O} \Big ( \L{p}{\nabla u_{\eps}}{\O} + \L{p}{\nabla g}{\O}\Big ) + \L{p}{\nabla u_{\eps}}{\O}^{p-1}\cdot\L{p}{g}{\O} \Bigg ]\\
&\le C(\mu,N,p, F, g,\O) \Bigg [ \L{p}{\nabla u_{\eps}}{\O} + \L{p}{\nabla u_{\eps}}{\O}^{p-1} +1\Bigg ].
\end{split}
\]
therefore, 
\be\label{uniform grad bound}
\L{p}{\nabla u_{\eps}}{\O}\le C(\mu,N,p, F, g,\O)
\ee
thus, we have the following uniform bound on $\L{p}{\nabla u_{\eps}}{\O}$.
{Moreover, from Poincaré inequality and \eqref{uniform grad bound} we have 
\be\label{uniform Lp bound}
\L{p}{u_{\eps}-g}{\O} \le C(\O) \L{p}{\nabla (u_{\eps}-g)}{\O}\implies \L{p}{u_{\eps}}{\O}  \le C(\mu,N,p, F, g,\O).
\ee}
From \eqref{uniform grad bound} and \eqref{uniform Lp bound}, we obtain that $u_{\eps}$ are uniformly bounded in $W^{1,p}(\O)$. Hence by Eberlyn's theorem, there exists at least one function $u\in W^{1,p}_g(\O)$ such that 
\be\label{weak limit}
u_{\eps}\wto u \;\;\mbox{weakly in $W^{1,p}(\O)$ as $\eps\to 0$}.
\ee
and due to compact embedding
\be\label{strong limit}
u_{\eps}\to u \mbox{ strongly in $L^p(\O)$ and } u_{\eps}\to u \mbox{ pointwise a.e. in $\O$}. 
\ee
\\
\textbf{Step 3. Proving that $u$ is a weak solution to \eqref{PDE}.}
\\
\\
Before delving into calculations, let us define some more notations:
\be\label{def2}
\begin{split}
U_{\eps^+} &:= \Psi^+_{\eps}(u_{\eps}), V_{\eps^+}:=\Psi^+_{\eps}(u),\\
U_{\eps^-} &:= \Psi^-_{\eps}(u_{\eps}),  V_{\eps^-} := \Psi^-_{\eps}(u).
\end{split}
\ee
We note that the following equalities clearly follows from definitions of $\Psi_{\eps}^{\pm}$
\be\label{split}
\begin{split}
u_{\eps} = U_{\eps^+} - U_{\eps^-}\\
\mbox{ and } u = V_{\eps^+} - V_{\eps^-}.
\end{split}
\ee 
Let us note two simple observations
\be\label{gradient of UV}
\begin{split}
\nabla U_{\eps^{\pm}} = \nabla \Big ( \Psi^{\pm}_{\eps}(u_{\eps}) \Big ) = \psi_{\eps}^{\pm}(u_{\eps})  \nabla u_{\eps}\\
\nabla V_{\eps^{\pm}} = \nabla \Big ( \Psi^{\pm}_{\eps}(u) \Big ) = \psi_{\eps}^{\pm}(u)  \nabla u
\end{split}
\ee
From the definition of $A_{\eps}(x,s)$ (c.f. \eqref{def1}) and from \eqref{gradient of UV} we observe that 
\be\label{def of Aeps}
\begin{split}
A_{\eps}(x,u_{\eps})|\nabla u_\eps|^{p-2}\nabla u_\eps &=\Big (A_+(x) ({\psi^+_{\eps}} (u_{\eps}) )^{p-1}+ A_-(x)({\psi^-_{\eps}} (u_{\eps}) )^{p-1} \Big ) |\nabla u_\eps|^{p-2}\nabla u_\eps\\ 
&= A_+(x) |\nabla U_{\eps^+}|^{p-2} \nabla U_{\eps^+}+A_-(x) |\nabla U_{\eps^-}|^{p-2} \nabla U_{\eps^-}
\end{split}
\ee

We claim that the following convergences are true up to a subsequence
\begin{enumerate}[label=\textbf{$($C\Roman*$)$}]
\item \label{convg 1}
$ U_{\eps^{\pm}} \wto u^{\pm} \mbox{ weakly in $W^{1,p}(\O)$} 
\mbox{ and } U_{\eps^{\pm}} \to u^{\pm} \mbox{ strongly in $L^p(\O)$ and a.e. in $\O$.}$

\item \label{convg 2}$ V_{\eps^{\pm}} \wto u^{\pm} \mbox{ weakly in $W^{1,p}(\O)$} 
\mbox{ and } V_{\eps^{\pm}} \to u^{\pm} \mbox{ strongly in $L^p(\O)$ and a.e. in $\O$.}
$
\item \label{convg 3}$ \big ({\psi_{\eps}^\pm} (u) \big )^{p-1}\nabla U_{\eps^\pm} \wto  \nabla u^{\pm}$ weakly in $L^p(\O)$.
\item \label{convg 4}$ \big ({\psi_{\eps}^\pm} (u) \big )^{p-1}\nabla V_{\eps^\pm} \wto  \nabla u^{\pm}$ weakly in $L^p(\O)$.
\item \label{convg 5}$ \big ({\psi_{\eps}^\pm} (u_\eps) \big )^{p-1}\nabla u_{\eps} \wto  \nabla u^{\pm}$ weakly in $L^p(\O)$.
\end{enumerate}
We look into the proofs of the claims \ref{convg 1}-\ref{convg 5} in the next step. For now, we assume the claims \ref{convg 1}-\ref{convg 5} to be true and use them to show that $u$ is a weak solution to \eqref{PDE} in $\O$. Firstly, we not that $u_\eps = u =g$ on $\partial \O$ in trace sense. Therefore, from the definitions of $U_{\eps^{\pm}}$ and $V_{\eps^{\pm}}$ (c.f. \eqref{def2}), it clearly follows that $U_{\eps^{\pm}} = V_{\eps^{\pm}}$ on $\partial \O$ in trace sense. In particular, $U_{\eps^+} -V_{\eps^+}\in W_0^{1,p}(\O)$. Since $u_\eps$ is a weak solution to \eqref{eps PDE}, we use $U_{\eps^+} -V_{\eps^+}\in W_0^{1,p}(\O)$ as a test function  
\be\label{weak def}
\underbrace{\int_{\O}\Big(  A_{\eps} (x,u_{\eps}) |\nabla u_{\eps}|^{p-2} \nabla u_{\eps} \cdot \nabla (U_{\eps^+} -V_{\eps^+})\Big )\,dx}_{\text{LHS}_\eps}  = \underbrace{\int_{\O} f_{\eps}(x,u_{\eps}) (U_{\eps^+} -V_{\eps^+})\,dx}_{\text{RHS}_\eps}
\ee
{With the same reasoning as in \eqref{controlling RHS 1} and \eqref{controlling RHS 2} and using \ref{convg 1} and \ref{convg 2}, we pass the limit $\eps\to 0$ on the right hand side of the equation above, we obtain 
$$
\int_{\O} f_{\eps}(x,u_{\eps})\cdot (U_{\eps^{+}} - V_{\eps^{+}}) \,dx \le C(N,p,\O) \L{N}{F}{\O} \L{p}{U_{\eps^{+}} - V_{\eps^{+}}}{\O}\to 0 \mbox{ as $\eps\to 0$.}
$$}
Therefore, 
\be\label{lim LHS}
\lim_{\eps\to 0} \text{LHS}_\eps  =\lim_{\eps\to 0}\int_{\O}\Big(  A_{\eps} (x,u_{\eps}) |\nabla u_{\eps}|^{p-2} \nabla u_{\eps} \cdot \nabla (U_{\eps^+} -V_{\eps^+})\Big )\,dx = 0
\ee
We consider the left hand side of the equation above. From  \eqref{split}, we have $(U_{\eps^+} -V_{\eps^+}) = (u_{\eps}-u) +(U_{\eps^-} -V_{\eps^-})$ and plugging it into the left side of \eqref{weak def} and using \eqref{def of Aeps}, we obtain
\[
\begin{split}
\text{LHS}_\eps &= \int_{\O}\Big(  A_{+} (x) |\nabla U_{\eps^+}|^{p-2} \nabla U_{\eps^+} \cdot \nabla (U_{\eps^+} -V_{\eps^+})\Big )\,dx + \int_{\O}\Big(  A_{-} (x) |\nabla U_{\eps^-}|^{p-2} \nabla U_{\eps^-} \cdot \nabla (U_{\eps^+} -V_{\eps^+})\Big )\,dx\\
&= \int_{\O}\Big(  A_{+} (x) |\nabla U_{\eps^+}|^{p-2} \nabla U_{\eps^+} \cdot \nabla (U_{\eps^+} -V_{\eps^+})\Big )\,dx + \int_{\O}\Big(  A_-(x) |\nabla U_{\eps^-}|^{p-2} \nabla U_{\eps^-} \cdot \nabla (U_{\eps^-}-V_{\eps^-})\Big )\,dx\\
& \qquad \qquad \qquad \qquad \qquad  \qquad \qquad \qquad \qquad \qquad \qquad+ \int_{\O}\Big(  A_- (x)|\nabla U_{\eps^-}|^{p-2} \nabla U_{\eps^-} \cdot \nabla (u_\eps-u )\,dx\\
&= \int_{\O}\Big(  A_{+} (x) |\nabla U_{\eps^+}|^{p-2} \nabla U_{\eps^+} \cdot \nabla (U_{\eps^+} -V_{\eps^+})\Big )\,dx + \int_{\O}\Big(  A_-(x) |\nabla U_{\eps^-}|^{p-2} \nabla U_{\eps^-} \cdot \nabla (U_{\eps^-}-V_{\eps^-})\Big )\,dx\\
& \qquad  \qquad \qquad  \qquad \qquad \qquad \qquad \qquad \qquad+ \int_{\O}\Big(  A_- (x)|\nabla u_{\eps}|^{p-2} \nabla u_{\eps} \cdot\Big ( (\psi_\eps^-(u_\eps))^{p-1} \nabla (u_\eps-u ) \Big )\,dx.\\
\end{split} 
\]
After doing a series of additions and subtractions in the above expression of $\text{LHS}_\eps$, we have 
\be\label{all integrals}
\begin{split}
\text{LHS}_\eps & =  \underbrace{\int_{\O} \Big ( A_+(x) (|\nabla U_{\eps^+}|^{p-2} \nabla U_{\eps^+} -  |\nabla V_{\eps^+}|^{p-2} \nabla V_{\eps^+})\cdot \nabla (U_{\eps^+} -V_{\eps^+})\Big )\,dx}_{\mathbf I 1 _{\eps}}\\
& \qquad \qquad \qquad \qquad \qquad \qquad \qquad \qquad \qquad + \underbrace{\int_{\O} \Big ( A_+(x) |\nabla V_{\eps^+}|^{p-2} \nabla V_{\eps^+}  \cdot \nabla (U_{\eps^+} -V_{\eps^+}) \Big )\,dx}_{\mathbf I2_\eps}\\
&\qquad + \underbrace{\int_{\O} \Big ( A_-(x) (|\nabla U_{\eps^-}|^{p-2} \nabla U_{\eps^-} -  |\nabla V_{\eps^-}|^{p-2} \nabla V_{\eps^-})\cdot \nabla (U_{\eps^-} -V_{\eps^-})\Big )\,dx}_{\mathbf {II}1_\eps}\\
& \qquad \qquad \qquad \qquad \qquad \qquad \qquad \qquad \qquad + \underbrace{ \int_{\O} \Big ( A_-(x) |\nabla V_{\eps^-}|^{p-2} \nabla V_{\eps^-}  \cdot \nabla (U_{\eps^-} -V_{\eps^-}) \Big )\,dx}_{\mathbf{II}2_\eps}\\
&+ \underbrace{\int_{\O}\Big(  A_- (x)\Big (  |\nabla u_{\eps}|^{p-2} \nabla u_{\eps} -|\nabla u|^{p-2} \nabla u \Big ) \cdot \Big( (\psi_\eps^-(u_\eps))^{p-1} \nabla (u_\eps-u ) \Big )\,dx}_{\mathbf{III}1_\eps}\\
&\qquad \qquad \qquad \qquad \qquad \qquad \qquad \qquad \qquad+ \underbrace{\int_{\O} A_-(x)|\nabla u|^{p-2} \nabla u  \cdot \Big( (\psi_\eps^-(u_\eps))^{p-1} \nabla (u_\eps-u ) \Big )\,dx}_{\mathbf{III}2_\eps} .
\end{split}
\ee
Let us look into $\mathbf{I}2_\eps$, $\mathbf{II}2_\eps$ and $\mathbf{III}2_\eps$. We claim that all three of those tend to $0$ as $\eps\to 0$. Indeed we have 
\[
\begin{split}
\mathbf{I}2_\eps &=\int_{\O} \Big ( A_+(x) |\nabla V_{\eps^+}|^{p-2} \nabla V_{\eps^+}  \cdot \nabla (U_{\eps^+} -V_{\eps^+}) \Big )\,dx \\
&= \int_{\O} \Big ( A_+(x) \Big ( |\nabla u|^{p-2} \nabla u \Big )  \cdot\Big ( \psi^+(u)^{p-1} \nabla (U_{\eps^+} -V_{\eps^+}) \Big ) \Big )\,dx
\end{split}
\]
We can easily check that since $\nabla u\in L^p(\O)$ and $0<\mu \le A_{\pm} \le \frac{1}{\mu}$, therefore $A_+(x) \Big ( |\nabla u|^{p-2} \nabla u \Big ) \in L^{p'}(\O)$. From \ref{convg 3} and \ref{convg 4}, we have $\psi^+(u)^{p-1} \nabla (U_{\eps^+} -V_{\eps^+}) \wto 0$ weakly in $L^p(\O)$. Hence,  we have
\be\label{I2}
\mathbf{I}2_\eps \to 0 \mbox{  as $\eps\to 0$ up to a subsequence}.
\ee
Because of the exact same reasoning as above we also have 
\be\label{II2}
\mathbf{II}2_\eps \to 0 \mbox{  as $\eps\to 0$ up to a subsequence}.
\ee 
Since we already know that $ A_-(x)|\nabla u|^{p-2} \nabla u \in L^{p'}(\O)$, therefore from \ref{convg 5}
\be\label{III2}
\mathbf{III}2_\eps \to 0 \mbox{  as $\eps\to 0$ up to a subsequence}.
\ee
Thus, using the information \eqref{lim LHS}, \eqref{I2}, \eqref{II2} and \eqref{III2} we rewrite \eqref{all integrals} as follows
\be
\lim_{\eps\to 0} \text{LHS}_\eps  = \lim_{\eps\to 0} (\mathbf{I}1_\eps +\mathbf{II}1_\eps +\mathbf{III}1_\eps ) =0
\ee
From convexity inequality $(|\xi|^{p-2}\xi - |\eta|^{p-2}\eta)\cdot(\xi-\eta)\ge 0$ for all $\xi,\eta\in \R^N$ and the fact that $A_{\pm}>\mu>0$, we have $\mathbf{I}1_\eps,\mathbf{II}1_\eps,\mathbf{III}1_\eps \ge 0$ for all $\eps>0$. Since addition of all three (non-negative) integrals $\mathbf{I}1_\eps,\mathbf{II}1_\eps$ and $\mathbf{III}1_\eps$ tend to zero, therefore all there integrals individually tend to zero. That is
\[
\begin{split}
\mathbf{I}1_\eps &\to 0\\
\mathbf{II}1_\eps &\to 0\\
\mathbf{III}1_\eps &\to 0
\end{split}
\]
In particular, we have 
\[
\begin{split}
\lim_{\eps\to 0} & \int_{\O} \Big ( A_+(x)  (|\nabla U_{\eps^+}|^{p-2} \nabla U_{\eps^+} -  |\nabla V_{\eps^+}|^{p-2} \nabla V_{\eps^+})\cdot \nabla (U_{\eps^+} -V_{\eps^+})\Big )\,dx =0 \mbox{ and }\\
\lim_{\eps\to 0} & \int_{\O} \Big ( A_-(x)  (|\nabla U_{\eps^-}|^{p-2} \nabla U_{\eps^-} -  |\nabla V_{\eps^-}|^{p-2} \nabla V_{\eps^-})\cdot \nabla (U_{\eps^-} -V_{\eps^-})\Big )\,dx=0.
\end{split}
\]
From the positivity of the integrand, this actually means
$$
A_\pm(x)  (|\nabla U_{\eps^\pm}|^{p-2} \nabla U_{\eps^\pm} -  |\nabla V_{\eps^\pm}|^{p-2} \nabla V_{\eps^\pm})\cdot \nabla (U_{\eps^\pm} -V_{\eps^\pm}) \to 0 \mbox{ in $L^1(\O)$. }
$$
Hence, up to a subsequence, we may assume that  
$$A_\pm(x)  (|\nabla U_{\eps^\pm}|^{p-2} \nabla U_{\eps^\pm} -  |\nabla V_{\eps^\pm}|^{p-2} \nabla V_{\eps^\pm})\cdot \nabla (U_{\eps^\pm} -V_{\eps^\pm}) \to 0\mbox{ a.e. in $\O$.}
$$
{Thus from \eqref{convex identity} and $\mu \le A_{\pm}\le \frac{1}{\mu}$ we have  (for $p\ge 2$)
$$
0\le \mu |\nabla (U_{\eps^\pm} - \nabla  V_{\eps^\pm})|^p \le  A_\pm(x)  (|\nabla U_{\eps^\pm}|^{p-2} \nabla U_{\eps^\pm} -  |\nabla V_{\eps^\pm}|^{p-2} \nabla V_{\eps^\pm})\cdot \nabla (U_{\eps^\pm} -V_{\eps^\pm}) \to 0.
$$
}
Therefore, we conclude that (up to a subsequence)
\be\label{pw convg}
\begin{split}
\nabla (U_{\eps^\pm} - \nabla  V_{\eps^\pm}) \to 0 \mbox{  a.e. as $\eps\to 0$}
\end{split}
\ee
Since $\psi_\eps^{\pm}\to \psi_0^{\pm}$ a.e. in $\R$, from the definition of $V_{\eps^{\pm}}$, $\nabla V_{\eps^{\pm}} = \psi_\eps^{\pm}(u) \nabla u \to \nabla u^{\pm}$  a.e. as $\eps\to 0$. Thus from \eqref{pw convg}
\be\label{the limit}
\begin{split}
\nabla U_{\eps^\pm}  \to \nabla u^{\pm} \mbox{  a.e. as $\eps\to 0$}
\end{split}
\ee
We already know that $|\nabla U_{\eps^{\pm}}|^{p-2}\nabla U_{\eps^{\pm}}$ is uniformly bounded in $L^{p'}(\O)$ (c.f. \eqref{uniform bound} in the next step) and therefore has a weak limit $\eps\to0$ in $L^{p'}(\O)$. This fact along with \eqref{the limit}, we obtain
$$
|\nabla U_{\eps^{\pm}}|^{p-2}\nabla U_{\eps^{\pm}} \wto |\nabla u^{\pm}|^{p-2}\nabla u^{\pm}\mbox{ weakly in $L^{p'}(\O)$}.
$$
Thus, for any $\vf\in C_c^{\infty}(\O)$, we have 
\be\label{LHS}
\begin{split}
\lim_{\eps\to 0}\int_{\O} \Big ( A_\pm(x) {(\psi^{\pm}}(u_{\eps}))^{p-1} |\nabla u_{\eps}|^{p-2}\nabla u_{\eps} \cdot \nabla \vf \Big )\,dx &= \lim_{\eps\to 0}\int_{\O} \Big ( A_\pm(x)   |\nabla U_{\eps^\pm}|^{p-2}\nabla U_{\eps^\pm} \cdot \nabla \vf \Big )\,dx\\
& =\int_{\O} \Big ( A_\pm(x)   {|\nabla u^\pm|}^{p-2}\nabla u^\pm \cdot \nabla \vf \Big )\,dx.
\end{split}
\ee 
Since $u_\eps \to u$ strongly in $L^p(\O)$, hence $u_{\eps}\to u$ pointwise a.e. up to subsequence. Thereofore ${\psi_\eps^{\pm} }(u_\eps)\to \psi_0^{\pm}(u)$ and therefore $f_{\eps}(x,u_{\eps}) \to f(x,u)$ pointwise a.e. in $\O$. Moreover, $|f_{\eps}(x,u_{\eps}) | \le |f_+ + f_-| \in L^1(\O)$ for all $\eps>0$, by Dominated convergence theorem, we have 
\be\label{RHS}
\lim_{\eps\to 0}\int_{\O} f_\eps(x,u_{\eps})\vf \,dx = \int_{\O}f(x,u)\vf\,dx.
\ee
We pass to the limit $\eps\to 0$ in \eqref{eps PDE weak} and using \eqref{LHS} and \eqref{RHS}, we obtain 
$$
\int_{\O} \Big ( A(x,u) {|\nabla u|}^{p-2}\nabla u \cdot \nabla \vf \Big )\,dx=\int_{\O}f(x,u)\vf\,dx.
$$
This shows, $u\in W_g^{1,p}(\O)$ is a weak solution to \eqref{PDE}. Now we finish the proof of Proposition \ref{existence} by proving the claims \ref{convg 1} -\ref{convg 5}.
\\
\\
\textbf{Step 4: Proving the claims \ref{convg 1} - \ref{convg 5}.}
\\
\\
Firstly, we observe that 
\be\label{auxiliary limit 0}
\begin{split}
\Psi_\eps^{\pm}(s) &\to s^{\pm} \mbox{ pointwise in $\R$.}\\
\psi_\eps^{\pm}(s) &\to \psi_0^{\pm}(s) \mbox{ pointwise in $\R$.}
\end{split}
\ee
Moreover, by compact embedding, $u_\eps \to u$ in $L^p(\O)$ therefore $u_\eps \to u$ almost everywhere in $\O$ up to a subsequence. Therefore we have 
\be\label{auxiliary limit 1}
U_{\eps^\pm} := \Psi_\eps^{\pm}(u_\eps) \to u^{\pm} \mbox{ up to a subsequence.}
\ee
Similarly, 
\be\label{auxiliary limit 2}
V_{\eps^\pm} := \Psi_\eps^{\pm}(u) \to u^{\pm} \mbox{ up to a subsequence.}
\ee 
On the other hand, we observe that
\be\label{auxiliary bounds}
\begin{split}
|\psi_\eps^{\pm}(s) | &\le 1,  \;\;\;\; \;\;\;\; \mbox{ for all $s\in \R$}\\
|\Psi_\eps^\pm(s)| &\le s^{\pm}+1, \mbox{ for all $s\in \R$.}
\end{split}
\ee
From  \eqref{uniform grad bound}, \eqref{uniform Lp bound} and \eqref{auxiliary bounds}, we obtain the following inequalities
\be\label{uniform bound}
\begin{split}
\int_{\O} |U_{\eps^\pm}|^p\,dx =  \int_{\O}|\Psi^\pm_{\eps}(u_{\eps})|^p\,dx &\le \int_{\O}|u_{\eps}^{\pm} +1|^p\,dx \le C(p) \int_{\O} |u_\eps +1|^p\,dx \le  C(N,\mu,p,F,g,\O).\\
\int_\O |\nabla U_{\eps^\pm}|^p\,dx =  \int_{\O}|\nabla \Psi^+_{\eps}(u_{\eps})|^p\,dx &\le \int_{\O}\psi^+_{\eps}(u_{\eps})^p |\nabla u_{\eps}|^p\,dx\le \int_{\O} |\nabla u_{\eps}|^p\,dx \le C(N,\mu,p,F,g,\O).\\
\int_{\O} |V_{\eps^\pm}|^p\,dx =  \int_{\O}|\Psi^\pm_{\eps}(u)|^p\,dx &\le \int_{\O}|u^{\pm} +1|^p\,dx \le C(p) \int_{\O} |u +1|^p\,dx \le  C(N,\mu,p,F,g,\O).\\
\int_\O |\nabla V_{\eps^\pm}|^p\,dx =  \int_{\O}|\nabla \Psi^+_{\eps}(u)|^p\,dx &\le \int_{\O}\psi^+_{\eps}(u)^p |\nabla u|^p\,dx\le \int_{\O} |\nabla u|^p\,dx \le C(N,\mu,p,F,g,\O).\\
\end{split}
\ee
The last two inequalities follows from the lower semicontinuity property of $W^{1,p}(\O)$ norm under weak convergence. This implies that $U_{\eps^{\pm}}$ and $V_{\eps^\pm}$ are uniformly bounded in $W^{1,p}(\O)$ and therefore, $U_{\eps^{\pm}}$ and $V_{\eps^\pm}$ have a weak limit in $W^{1,p}(\O)$ up to a subsequence. From \eqref{auxiliary limit 1} and \eqref{auxiliary limit 2} we conclude the following 
\be\label{auxiliary weak limit}
\begin{split}
U_{\eps^\pm} &\wto u^\pm \mbox{ weakly in $W^{1,p}(\O)$}\\
V_{\eps^\pm} &\wto u^\pm \mbox{ weakly in $W^{1,p}(\O)$}\\
\end{split}
\ee
Thus, we prove claims \ref{convg 1} and \ref{convg 2}. In order to prove \ref{convg 3} and \ref{convg 4}, first we observe that from \eqref{auxiliary limit 0} $(\psi_{\eps}^\pm(u))^{p-1} \to (\psi_0^\pm(u))^{p-1}$ almost everywhere in $\O$. We pick any $\bf \in L^p(\O)^N$ and therefore we have 
\be\label{auxiliary convergence 3}
(\psi_{\eps}^\pm(u))^{p-1} \bf \to (\psi_0^\pm(u))^{p-1} \bf \mbox{ almost everywhere in $\O$.}
\ee
Since $|(\psi_{\eps}^\pm(u))^{p-1} \bf| \le |\bf|\in L^p(\O)$, by Dominated Convergence Theorem (DCT) and \eqref{auxiliary convergence 3} we conclude 
\be\label{auxiliary convergence 4}
(\psi_{\eps}^\pm(u))^{p-1} \bf \to (\psi_0^\pm(u))^{p-1} \bf \mbox{ in $L^p(\O)$.}
\ee
From \ref{convg 1} and \cite[Proposition 3.13 (iv)]{HB10} we have
\[
\begin{split}
\int_{\O} \nabla U_{\eps^\pm}\cdot \Big ( (\psi_{\eps}^\pm(u))^{p-1} \bf \Big )\,dx  \to \int_\O \nabla u^\pm \cdot \Big ( (\psi_{0}^\pm(u))^{p-1} \bf \Big )
\end{split}
\]
in other words 
\be\label{auxiliary limit 5}
\int_{\O} \Big ( (\psi_{\eps}^\pm(u))^{p-1} \nabla U_{\eps^\pm} \Big )\cdot \bf\,dx \to \int_{\O} \nabla u^{\pm}\cdot \bf\,dx \mbox{ for all $\bf\in L^p(\O)^N$}.
\ee
This proves \ref{convg 3}. The proof of \ref{convg 4} follows via the exact same reasoning as above using \eqref{auxiliary convergence 4}, \ref{convg 2} and \cite[Proposition 3.13 (iv)]{HB10}. In order to prove \ref{convg 5}, we observe that from \eqref{auxiliary limit 0} and the fact $u_\eps \to u$ almost everywhere in $\O$, $(\psi_\eps^{\pm}(u_\eps))^{p-1}\bf \to (\psi_0^\pm(u))^{p-1} \bf$ almost everywhere in $\O$ for all $\bf \in L^p(\O)^N$, hence by DCT $(\psi_\eps^{\pm}(u_\eps))^{p-1}\bf \to (\psi_0^\pm(u))^{p-1} \bf$ in $L^p(\O)^N$. We use the fact that $\nabla u_\eps \wto \nabla u$ weakly in $L^p(\O)$, the proof of \ref{convg 5} follows from the same reasoning as in the proof of \ref{convg 3} and \ref{convg 4}. 
\end{proof}
}

\section{Approximation lemma via compactness}\label{compactness}
\begin{remark}\label{rescaling 0.1}
For a given function $u$ satisfying the PDE \eqref{PDE} in $B_{\Theta}(x_0)$ we define $w$ as follows
$$
w(y) := \Phi u(\Theta y +x_0) + \Psi,  \;\;y\in B_1.
$$
By rescaling and change of variables, we can easily verify that the new rescaled function $w$ satisfies the following PDE in $B_1$
$$
-\dive(\bar A(x,w)|\nabla w|^{p-2} \nabla w) = \bar f(x,w) 
$$
where $\bar A_{\pm}$ and $\bar f_{\pm}$ are defined as follows
\[
\begin{split}
\bar A_{\pm}(x) &= A(\Theta x +x_0)\\
\bar f_{\pm}(x)&= \Phi^{p-1} \Theta ^{p} f(\Theta x +x_0).
\end{split}
\]
\end{remark}
\begin{proposition}\label{step 1}
Suppose $u\in W^{1,p}(B_{1/2})$ is a weak solution of \eqref{PDE} in $B_{1/2}$ such that we have $\Linfty{u}{B_{1/2}} \le 1$ and $\L{p}{\nabla u}{B_{1/2}}\le M$. Then for every $\eps>0$ there exists $\delta (\eps, N,p, \mu, M) >0$ such that if 
$$
\max \Big ( \Linfty{A_{\pm} - A_{\pm}(0)}{B_{1/2}}, \LN{f_{\pm}}{B_{1/2}} \Big ) \le \delta
$$
then 
$$
\Linfty{u-h}{B_{1/4}}\le \eps
$$
for some $h\in W^{1,p}(B_{1/2})$ such that 
\be\label{limit PDE}
\dive (A(0, h) |\nabla h|^{p-2}\nabla h) =0 \;\mbox{in $B_{1/2}$}.
\ee
\end{proposition}

\begin{proof}
Let us suppose by contradiction that the statement of Proposition \ref{step 1} is not true. This implies that there exists $\eps_0>0$ and a sequence $A_k, f_k$ such that $A_{\pm,k}\in C(B_{1/2})$ and $f_k\in L^{N}(B_{1/2})$ and $\|A_{\pm,k}-A_{\pm}(0)\|_{L^{\infty}(B_{1/2})}<\frac{1}{k}$, $\|f_{\pm,k}\|_{L^N(B_{1/2})}<\frac{1}{k}$. As well as for corresponding weak solutions $u_k$ of \eqref{PDEk} (see below) such that $\L{p}{\nabla u_k}{B_{1/2}}\le M$ and $\| u_k \|_{L^{\infty}(B_{1/2})} \le 1$ and for every $h$ satisfying \eqref{limit PDE} we have
\be\label{absurdum}
\|u_k-h\|_{L^{\infty}(B_{1/4})}>\eps_0.
\ee
$u_k$ are weak solution to the following PDE
\be\label{PDEk}
\dive(A_k(x,u_k)|\nabla u_k|^{p-2}\nabla u_k) = f_k(x,u_k),\;\;\mbox{in $B_1$}.
\ee
We know that the PDEs \eqref{PDEk} satisfy the structural condition in \cite[Chapter 10, Section 1]{ED09} and therefore $u_k$ belong to De-Giorgi class $DG_p(\mu, N) $. Therefore, $u_k$ are locally bounded in $B_1$ (c.f. \cite[Chapter 10, Theorem 2.1]{ED09}) and also Hölder continuous in $B_{1/4}$( c.f. \cite[Chapter 10, Theorem 3.1]{ED09}). 
It follows the existence of $\beta_0:= \beta_0 ( p,\mu,N)$ and $C_0:= C_0 (p,\mu,N)$ such that for every $k\in \N$ we have
$$
\| u_k\|_{C^{\beta_0}(B_{1/4})}  \le C_0.
$$
By Arzela Ascoli theorem, there exists $u_0\in C^{0,\alpha_0}(B_{1/4})$ such that 
\be\label{uniform convergence p=2} 
u_k\to u_0 \mbox{ in $L^{\infty}(B_{1/4})$ up to a subsequence. }
\ee
Since $\L{p}{\nabla u_k}{B_{1/2}}\le M$ and $\L{p}{u_k}{B_{1/2}} \le C(N,p)\cdot \Linfty{u_k}{B_{1/2}}\le C(N,p)$ for every $k$, $u_k$ is a bounded sequence in $W^{1,p}(B_{1/2})$.  That is 
$$
\| u_k \|_{W^{1,p}(B_{1/2})} \le C(M, N,p).
$$
Therefore, $u_k$ converges weakly to $u_0$ in $W^{1,p}(B_{1/4})$ (up to a subsequence)
$$
u_k \wto u_0 \;\;\mbox{weakly in $W^{1,p}(B_{1/4})$}.
$$
Since $u_k$ are weak solutions of \eqref{PDEk}, for every $\Phi \in W_0^{1,p}(B_{1/2})$ we have
\be\label{E1}
\int_{B_{1/2}} \big ( A_k(x,u_k) |\nabla u_k|^{p-2}\nabla u_k\cdot \nabla \Phi \big )\,dx = \int_{B_{1/2}}f_k \Phi \,dx
\ee 
and therefore we have 
\be\label{E2}
\int_{B_{1/2}} \Big (  \big ( A_k(x,u_k)-A(0,u_k) \big ) |\nabla u_k|^{p-2}\nabla u_k\cdot \nabla \Phi \Big )\,dx + \int_{B_{1/2}}A(0,u_k ) |\nabla u_k|^{p-2}\nabla u_k\cdot \nabla \Phi \,dx  = \int_{B_{1/2}}f_k \Phi \,dx
\ee
Sinnce $A_{k, \pm} \to A_{\pm}(0)$ uniformly in $B_{1/2}$ we have 
\[
\begin{split}
\int_{B_{1/2}} \Big (  \big ( A_k(x,u_k)-A(0,u_k) \big ) |\nabla u_k|^{p-2}\nabla u_k\cdot \nabla \Phi \Big )\,dx & \le \Linfty{A_k-A(0)}{B_{1/2}} \int_{B_{1/2}}|\nabla u_k|^{p-2}\nabla u_k\cdot \nabla \Phi \,dx\\
&\le \frac{1}{k} \L{p}{\nabla u_k}{B_{1/2}}^{p-1}\L{p}{\nabla \Phi}{B_{1/2}} \\
&\le \frac{1}{k}M^{p-1}\L{p}{\nabla \Phi}{B_{1/2}} \to 0   \;\;\mbox{as $k\to \infty$}.
\end{split}
\]
Moreover we also observe that if we define $F_k:= |f_{+,k}| + |f_{-,k}|$. For $N>p$ we have ${p^*}' =\frac{Np'}{N+p'}<N$ and by Sobolev embedding we have
\[
\begin{split}
\int_{B_{1/2}} |f_k(x,u_k)\Phi |\,dx &\le \int_{B_{1/2}} (|f_{+,k}| + |f_{-,k}| )|\Phi | \,dx = \int_{B_{1/2}} |F_k| |\Phi |\,dx \\
&\le \|F_k\|_{L^{{p^*}'}(B_{1/2})} \|\Phi \|_{L^{p^*}(B_{1/2})}\\
&\le  \frac{C(p,N)}{k}\|\Phi \|_{W^{1,p}(B_{1/2})}\to 0.\\
\end{split}
\]
For $p\ge N$, we know that $\Phi \in W_0^{1,p}(B_{1/2}) \subset L^{q}(B_{1/2})$ for every $q\ge 1$. In particular $\Phi \in L^{N'}(B_{1/2})$. Therefore
\[
\begin{split}
\int_{B_{1/2}} |f_k(x,u_k)\Phi |\,dx &\le \int_{B_{1/2}} (|f_{+,k}| + |f_{-,k}| )|\Phi | \,dx = \int_{B_{1/2}} |F_k| |\Phi |\,dx \\
&\le C(N) \|F_k\|_{L^{N}(B_{1/2})} \|\Phi \|_{L^{N'}(B_{1/2})}\\
&\le  \frac{C(N,p)}{k}\|\Phi \|_{W^{1,p}(B_{1/2})}\to 0.\\
\end{split}
\]
Plugging the above computations in \eqref{E2} we obtain 
\be\label{E3}
\lim_{k\to \infty} \int_{B_{1/2}}A(0,u_k ) |\nabla u_k|^{p-2}\nabla u_k\cdot \nabla \Phi \,dx =0.
\ee
\textbf{We claim that $u_0$ satisfy the PDE \eqref{limit PDE} which will give us a contradiction and prove Proposition \ref{step 1}.} We ease the notation in upcoming computations by renaming $A_{\pm}(0)$ as follows
\[
\begin{split}
A_+(0)&=: a^{p-1}\\
A_-(0)& =: b^{p-1}
\end{split}
\]

Let $\eta\in C_c^{\infty}(B_{1/2})$ be such that 
$$
\eta= \begin{cases}
1\;\;\mbox{ in $B_{\frac{1}{4}}$}\\
0\;\;\mbox{ on $\partial B_{1/2}$}.
\end{cases}
$$
Let us consider $\bf_k \in W_0^{1,p}(B_{1/2})$ such that we have
\be\label{test function}
\bf_k := \eta  (\TT_{a,b}(u_k) - \TT_{a,b}(u_0)). 
\ee
We also rename the following functions to ease computations below
\[
\begin{split}
U_k &:= \TT_{a,b}(u_k)\\
U_0 &:= \TT_{a,b}(u_0).
\end{split}
\]
We can verify that in order to prove that $u_0$ satisfies the PDE \eqref{limit PDE} it is enough to show that $U_0$ is $p$-harmonic in $B_{1/2}$. Indeed from \eqref{gradient relation}, we have 
\[
\begin{split}
\dive(A(0,u_0)|\nabla u_0|^{p-2}\nabla u_0)  &=\dive(a^{p-1}|\nabla u_0^+|^{p-2}\nabla u_0^+ -b^{p-1}|\nabla u_0^-|^{p-2}\nabla u_0^-)
\\ &= \dive(|\nabla (\TT_{a,b}(u))|^{p-2} \nabla (\TT_{a,b}(u)))=\Delta_p U_0.
\end{split}\]
This means, $u_0$ satisfying \eqref{limit PDE} is equivalent to
$$
\Delta_p U_0 = 0 \;\;\mbox{in $B_{1/2}$}. 
$$
From \eqref{test function}, we can write 
\be\label{test function 2}
\Phi_k = \eta( U_k-U_0).
\ee
Since the map $v\to \TT_{a,b}(v)$ is sequentially continuous in strong and weak $W^{1,p}(B_{1/2})$ topology \cite[Proposition 3.7 (1c), Proposition 3.7 (2c)]{MS22}. In short, the following convergences hold up to a subsequence
\begin{enumerate}[label = \textbf{(C\arabic*)}, ref= C\arabic*]
\item \label{C1} $u_k\wto u_0$ and $U_k \wto U_0$ weakly in $W^{1,p}(B_{1/2})$.
\item \label{C2} $u_k \to u_0$ and $U_k \to U_0$ strongly in $L^p(B_{1/2})$.
\item \label{C3} $u_k\to u_0$ and $U_k \to U_0$ pointwise almost everywhere in $B_{1/2}$.
\end{enumerate}
We claim that $\nabla U_k \to \nabla U_0$ alomost everywhere in $B_{1/2}$. In order to prove it, we observe that $A(0,u_k)|\nabla u_k|^{p-2} \nabla u_k = |\nabla (\TT_{a,b}(u_k))|^{p-2}\nabla (\TT_{a,b}(u_k))=|\nabla U_k|^{p-2} \nabla U_k$ (c.f. \eqref{gradient relation}). Therefore, we can write \eqref{E3} 
\be\label{e4}
\lim_{k\to \infty} \int_{B_{1/2}} |\nabla U_k|^{p-2} \nabla U_k \cdot \nabla \Phi \,dx =0.
\ee
for every $\Phi\in W_0^{1,p}(B_{1/2})$. Now, we take $\Phi: = \Phi_k$  in the previous identity \eqref{e4} and plug  in the definition of $\Phi_k$ from  \eqref{test function 2}
\be\label{e5}
\lim_{k\to \infty} \Big [ \int_{B_{1/2}} (U_k-U_0) |\nabla U_k|^{p-2}\nabla U_k\cdot \nabla \eta \,dx + \int_{B_{1/2}} \eta |\nabla U_k|^{p-2}\nabla U_k \cdot \nabla (U_k-U_0)\,dx  \Big ]=0.
\ee
We observe that the sequence $|\nabla U_k|^{p-2}\nabla U_k\cdot \nabla \eta$ is bounded in $L^{p'}(B_{1/2})$ and from \eqref{C2} we have 
$$
\lim_{k\to \infty}\int_{B_{1/2}} (U_k-U_0) |\nabla U_k|^{p-2}\nabla U_k\cdot \nabla \eta \,dx  =0.
$$
and using the above identity in \eqref{e5} we obtain 
\be\label{e6}
\lim_{k\to \infty}\int_{B_{1/2}} \eta |\nabla U_k|^{p-2}\nabla U_k \cdot \nabla (U_k-U_0)\,dx  =0.
\ee
We rewrite \eqref{e6} as follows
\be\label{e7}
\lim_{k \to \infty} \Bigg [ \int_{B_{1/2}} \Big (  \eta  \big ( |\nabla U_k|^{p-2}\nabla U_k  - |\nabla U_0|^{p-2} \nabla U_0\big ) \cdot \nabla (U_k-U_0) \Big ) \,dx + \int_{B_{1/2}} \eta |\nabla U_0|^{p-2} \nabla U_0\cdot \nabla (U_k-U_0)\,dx \Bigg ] =0
\ee
Since $\eta |\nabla U_0|^{p-2} \nabla U_0 \in L^{p'}(B_{1/2})$, from \eqref{C1} we have 
$$
\lim_{k \to \infty} \int_{B_{1/2}} \eta |\nabla U_0|^{p-2} \nabla U_0\cdot \nabla (U_k-U_0)\,dx =0
$$
and therefore from \eqref{e7}
\be\label{e8}
\lim_{k \to \infty}  \int_{B_{1/2}} \eta \Big (   \big ( |\nabla U_k|^{p-2}\nabla U_k  - |\nabla U_0|^{p-2} \nabla U_0\big ) \cdot \nabla (U_k-U_0) \Big ) \,dx =0.
\ee

For ease of notation, let us define $\GG_k$ as follows
$$
\GG_{k}(x):= \eta  \big ( |\nabla U_k|^{p-2}\nabla U_k  - |\nabla U_0|^{p-2} \nabla U_0\big ) \cdot \nabla (U_k-U_0)
$$
From \eqref{convex identity}, $\GG_k\ge 0$ and from \eqref{e8}, we have $\GG_k \to 0$ in $L^1(B_{1/2})$. Therefore, we conclude the following up to a subsequence
\be\label{e9}
\GG_{k} \to 0 \;\mbox{ pointwise a.e.  in $B_{1/2}$.}
\ee 
From Lemma \ref{claim} we have 
\be
\mbox{$\nabla U_k \to \nabla U_0$ pointwise almost everywhere in $B_{1/2}$. }
\ee
 From this, we deduce that 
\begin{enumerate}[label = {(\textbf{R\arabic*})}]
\item \label{R1} $|\nabla U_k|^{p-2} \nabla U_k \to |\nabla U_0|^{p-2} \nabla U_0$ pointwise a.e. in $B_{1/2}$ up to a subsequence.
\item \label{R2} Since $|\nabla U_k|^{p-2} \nabla U_k$ is bounded sequence in $L^{p'}(B_{1/2})$ and therefore from \ref{R1} and {Eberlyn's theorem} we have
$$
|\nabla U_k|^{p-2} \nabla U_k \wto |\nabla U_0|^{p-2} \nabla U_0\mbox{ weakly in $L^{p'}(B_{1/2})$.}
$$
\end{enumerate}
Now, we again look at \eqref{e4}. For $\Phi\in W_0^{1,p}(B_{1/2})$ we have 
$$
\lim_{k\to \infty} \Big ( \int_{B_{1/2}} \big ( |\nabla U_k|^{p-2}\nabla U_k -|\nabla U_0|^{p-2} \cdot \nabla U_0 \big )\cdot \nabla \Phi \,dx +  \int_{B_{1/2}}|\nabla U_0|^{p-2} \cdot \nabla U_0 \cdot \nabla \Phi \,dx \Big )  =0.
$$

From \ref{R2} the first integral in the above limit tends to zero as $k$ tends to $\infty$. Therefore
$$
\int_{B_{1/2}}|\nabla U_0|^{p-2} \cdot \nabla U_0 \cdot \nabla \Phi \,dx =0, \;\forall\, \Phi \in W_0^{1,p}(B_{1/2}).
$$
This proves that $U_0$ is $p$-harmonic in $B_{1/2}$ and equivalently $u_0$ satisfies the PDE \eqref{limit PDE}.
\end{proof}
{
Now, we remove the dependence of $\delta(\eps, N,p, \mu, M)$ in Proposition \ref{step 1} from $M$ (the bound on $\int_{B_{1/2}}|\nabla u|^p\,dx$). We  do so via Widman's hole filling technique.
}

\begin{proposition}\label{step 1.5}
Suppose $u\in W^{1,p}(B_{1})$ is a weak solution of \eqref{PDE} in $B_1$ such that we have $\Linfty{u}{B_{1}} \le 1$. Then for every $\eps>0$ there exists $0<\delta (\eps, N,p, \mu) <1$ such that if 
$$
\max \Big ( \Linfty{A_{\pm} - A_{\pm}(0)}{B_{1}}, \LN{f_{\pm}}{B_{1}} \Big ) \le \delta
$$
then 
$$
\Linfty{u-h}{B_{1/4}}\le \eps
$$
for some $h\in W^{1,p}(B_{1/4})$ satisfying \eqref{limit PDE}.
\end{proposition}

\begin{proof}
We claim that $\L{p}{\nabla u}{B_{1/2}} \le M(N,p, \mu)$. In order to prove this bound, let $r,s>0$ be such that $1/2\le s<t \le 1$ and $\eta \in C_0^{\infty}(B_1)$ such that $0\le \eta \le 1$ and
$$
\eta(x)=
\begin{cases}
1 \;\;x\in B_s\\
0\;\;x\in B_1\setminus B_t.
\end{cases}
$$
we can assume that 
\be\label{grad eta}
|\nabla \eta| \le \frac{C(N)}{|s-t|}
\ee
We consider $\vf = \eta u \in W_0^{1,p}(B_1)$ and from \eqref{weak solution} ($\O = B_1$) we have
\[
\begin{split}
\int_{B_1} \Big ( A(x,u)|\nabla u|^{p-2}\nabla u\cdot \nabla (\eta u) \Big )\,dx = \int_{B_1} f(x,u) (\eta u)\,dx.
\end{split}
\]
We expand the integral on LHS and we obtain 
\be\label{f1}
\int_{B_t} \eta A(x,u) |\nabla u|^p \,dx + \int_{B_t} \big ( u A(x,u) |\nabla u|^{p-2}\nabla u\cdot \nabla \eta  \big )\,dx  = \int_{B_t} f(x,u) (\eta u)\,dx.
\ee
Since $\nabla \eta$ is supported in $B_t\setminus B_s$,
\be\label{f1.1}
\begin{split}
\int_{B_t} \big ( u A(x,u) |\nabla u|^{p-2}\nabla u\cdot \nabla \eta  \big )\,dx &=\int_{B_t\setminus B_s} \big ( u A(x,u) |\nabla u|^{p-2}\nabla u\cdot \nabla \eta  \big )\,dx\\
&\le \frac{1}{\mu} \L{p}{\nabla u}{B_t\setminus B_s}^{p-1} \L{p}{\nabla \eta }{B_1}
\end{split}
\ee
we define $F := |f_+| + |f_-|$ and we have 
\be\label{f1.2}
\begin{split}
\int_{B_1} \Big | f(x,u) (\eta u)\,dx \Big |&= \int_{B_t} f(x,u) (\eta u)\,dx\\
& \le  C(N) \LN{F}{B_1} \le C(N).
\end{split}
\ee
Since $\eta A(x,u) |\nabla u|^p\ge 0$ (from \ref{H3}), \eqref{f1} leads to the following inequality
$$
\int_{B_s} A(x,u) |\nabla u|^p  \,dx\le \int_{B_t}  \eta A(x,u) |\nabla u|^p  \,dx \le  \int_{B_t} \big | \big  ( u A(x,u) |\nabla u|^{p-2}\nabla u\cdot \nabla \eta  \big ) \big | \,dx  + \int_{B_t} \big | f(x,u) \big | \big | (\eta u) \big |\,dx.
$$
Combining \eqref{f1.1}, \eqref{f1.2} and the ellipticity assumption \ref{H3} with the inequality above, for any $\delta_0>0$ we have 
\[
\begin{split}
\mu \int_{B_s} |\nabla u|^p\,dx &\le \frac{1}{\mu} \L{p}{\nabla u}{B_t\setminus B_s}^{p-1} \L{p}{\nabla \eta }{p} + \Big | \int_{B_t} f(x,u) (\eta u)\,dx \Big |\\
&\le \frac{C(\mu, p)}{ {\delta_0} ^{p'}} \L{p}{\nabla u}{B_t\setminus B_s}^{p} + {\delta_0} ^{p} \L{p}{\nabla \eta }{p} ^p +C(N).
\end{split}
\]
By taking $\delta_0=1/2$, from \eqref{grad eta} we obtain 
\be\label{f2}
\begin{split}
\int_{B_s} |\nabla u|^p\,dx &\le C_1(p,\mu)\int_{B_t \setminus B_s} |\nabla u|^p\,dx  + \frac{C_2(N)}{|s-t|^p} +C(N).
\end{split}
\ee
We add the term $C_1\int_{B_s} |\nabla u|^p\,dx$ on both sides of \eqref{f2} and we arrive at 
\be\label{f3}
\begin{split}
\int_{B_s} |\nabla u|^p\,dx \le \frac{C_1}{C_1+1}\int_{B_t} |\nabla u|^p\,dx + \frac{C_2}{|s-t|^p}+ C(N).
\end{split}
\ee
Now, from \cite[Lemma 6.1]{eg05}, we have
\be\label{f4}
\int_{B_{1/2}} |\nabla u|^p\,dx \le C_3(N,p, \mu ).
\ee
Now we apply Proposition \ref{step 1}. Since $u$ is a weak solution to \eqref{PDE} in $B_{1}$ and hence in ${B_{1/2}}$.  We can choose $\delta(\eps,  M, \mu, p, N) >0$ in Proposition \ref{step 1} which corresponds to $M= C_3(p,\mu,N)^{1/p}$. Therefore we have $\delta:= \delta(\eps, p,\mu,N)$ such that Proposition \ref{step 1.5} holds.
\end{proof}
\begin{remark}
From \eqref{gradient relation} and \eqref{limit PDE}, we know that $\Delta_p(\TT_{a,b}(h)) =0$ in $B_{1/2}$ ($h$ as in Lemma \ref{step 1.5}). Since $p$-harmonic functions are locally $C^{1,\gamma}$ regular for {some $\gamma:=\gamma(p,N)$.} From Lemma \ref{preservation of regularity} we write 
\be\label{limit regularity}
h\in C^{0,1}(B_{1/4}).
\ee
\end{remark}
\section{Optimal regularity of weak solutions}\label{regularity}
In this section, we follow the same steps as in the proofs of regularity theory for minimizers in \cite[Section 7]{MS22}, we adapt those proofs in the context of weak solutions to \eqref{PDE} and $p\in (1,\infty)$.
\begin{lemma}\label{step 2}
Suppose $u\in W^{1,p}(B_1)$ weakly solves the PDE \eqref{PDE} in $B_1$ with $\| u\|_{L^{\infty}(B_1)}\le 1$ and $u(0)=0$.  Then for any $0<\alpha<1$, there exists $\delta (N,p,\mu,\alpha) >0$ and $0< R_0(N,p,\mu,\alpha)<1/4$ such that if 
$$
 \max \Big (\|A_{\pm}-A_{\pm}(0)\|_{L^{\infty}(B_1)},\, \|f_{\pm}\|_{L^N(B_1)} \Big ) < \delta 
$$
then we have 
\be\label{claim step 2}
\sup_{B_{R_0}}|u-u(0)|\le R_0^{\alpha}.
\ee
\end{lemma}
\begin{proof}
Let $\eps>0$ which will be suitably chosen later. We know that for $\delta(\eps, N,p, \mu)>0$ and $h\in W^{1,p}(B_{1/4})$ as in Proposition \ref{step 1.5} we have  
\be\label{s2e1}
\|u-h\|_{L^{\infty}(B_{1/4})}<\eps.
\ee 
Fix $\beta=\frac{1+\alpha}{2}$ from \eqref{limit regularity} we have 
\be\label{s2e2}
\sup_{B_r}|h-h(0)| \le C(N,p,\mu, \alpha)r^{\beta}\;\;\forall r<1/4.
\ee
Combining equations (\ref{s2e1}) and (\ref{s2e2}) we get for $r<1/4$
\be\label{s2e3}
\begin{split}
\sup_{B_r}|u(x)-u(0)|&\le \sup_{B_r}\Big ( |u(x)-h(x)|+|h(x)-h(0)|+|h(0)-u(0)| \Big ) \\
&\le 2\eps+C(N,p,\mu, \alpha)r^{\beta}.
\end{split}
\ee
In the equation above, we select $r=R_0(N,p,\mu,\alpha)<1/4$ such that 
$$
C(N,p, \mu,\alpha)R_0^{\beta}=\frac{R_0^{\alpha}}{3}
$$
that is 
$$
R_0=\Big (  \frac{1}{3C} \Big )^{2/(1-\alpha)}.
$$
Now, we choose $\eps(N,p,\mu, \alpha)$ in such a way that 
$$
\eps<\frac{R_0^{\alpha}}{3}.
$$
We see that the choice of $\delta$ depending on $\eps$ and since $\eps$ depends on $N,p,\mu$ and $\alpha$ therefore $\delta$ is actually chosen depending on $N,p,\mu$ and $\alpha$. We use the fact that $u(0)=0$ and since $C(N,p,\mu, \alpha)R_0^{\beta}$ and $\eps$ are bounded by $R_0^{\alpha}/3$. From \eqref{s2e3} we have 
$$
\sup_{B_{R_0}} |u| \le R_0^{\alpha}.
$$
\end{proof} 
We now have ingredients to show the $C^{0,1^-}$ estimates for a minimizer $u$ around the set $\{u=0\}$, in particular, the free boundary $F(u)$.

\begin{lemma}\label{step 3}
Suppose $u$ satisfy the PDE \eqref{PDE}  with $\| u\|_{L^{\infty}(B_1)}\le 1$  and $u(0)=0$.  Then for all $0<\alpha<1$ and $\delta(N,p, \mu,\alpha)>0$ as in Lemma \ref{step 2} there exists $C(N,p,\mu,\alpha)>0$ such that if 
$$
 \max \Big (\|A_{\pm}-A_{\pm}(0)\|_{L^{\infty}(B_1)},\, \|f_{\pm}\|_{L^N(B_1)} \Big ) < \delta 
$$
then for $R_0(N,p, \alpha)$ as in Lemma \ref{step 2} we have
\be\label{claim step 3.}
\sup_{B_{r}}|u(x)|\le C(N,p, \mu,\alpha) \cdot r^{\alpha} \;\; \forall r<R_0.
\ee
Precisely speaking, we have $C(N,p,\mu,\alpha)= R_0^{-\alpha}$.
\end{lemma}

\begin{proof}
We argue by induction and rescaling, we claim that 
\be\label{s3e1}
 \sup_{B_{R_0^{k}}}|u(x)|\le R_0^{k\alpha }\;\;\;\forall k\in \N.
\ee
From Lemma \ref{step 2} we can see that (\ref{s3e1}) holds for $k=1$, and suppose it holds up to $k_0\in \N$. We prove that (\ref{s3e1}) holds for $k=k_0+1$. We define the following rescaling 
$$
\tilde u(y)=\frac{1}{R_0^{k_0\alpha}}u(R_0^{k_0}y).
$$
From Remark \ref{rescaling 0.1} we have
\be\label{rescaled PDE}
\dive(\tilde A(x,\tilde u)|\nabla \tilde u|^{p-2} \nabla \tilde u) = \tilde f(x,\tilde u) \;\;\mbox{in $B_1$}
\ee
with
\[
\begin{split}
&\tilde A_{\pm}(y,s)=A_{\pm}(R_0^{k_0}y,s)\\
&\tilde f_{\pm}(y,s)=R_0^{(pk_0(1-\alpha)+k_0\alpha)}f_{\pm}(R_0^{k_0}y,s).\\
\end{split}
\]
We verify that the functional $\tilde u$ satisfies the assumptions of Lemma \ref{step 2}. Indeed, from \eqref{s3e1}, we have 
$$
\sup _{B_{1}}|\tilde u|= R_0^{-k_0 \alpha} \sup_{B_{R_0^{k_0}}}|u|\le 1
$$
also for $\delta$ as in Lemma \ref{step 2} we can see that
$$
\sup_{B_1}|\tilde A_{\pm} - \tilde A_{\pm}(0)|=  \sup_{B_{R_0^{k_0}}} |A_{\pm}  -A_{\pm}(0)| \le \delta
$$
and 
$$
\|\tilde f_{\pm}\|_{L^{N}(B_1)}=R_0^{k_0(1-\alpha)(p-1)}\|f_{\pm}\|_{L^N(B_{R_0^k})}\le  \delta.
$$
Moreover $\tilde u(0)=0$, and hence we verify all the assumptions of Lemma \ref{step 2} for $\tilde u$. Therefore,
$$
\sup_{B_{R_0}}|\tilde u|\le R_0^{\alpha}
$$
on putting back the definition of $\tilde{u}$, we obtain the equation above in terms of $u$ 
$$
\sup_{B_{R_0^{k_0+1}}}|u|\le R_0^{(k_0+1)\alpha}.
$$
Hence we have proven the claim (\ref{s3e1}). To prove (\ref{claim step 3.}), we fix $0<r<R_0$ and $k(r)\in \N$ such that $R_0^{k+1}\le r<R_0^k$. From (\ref{s3e1}) we see that 
$$
\sup_{B_r}|u|\le \sup_{B_{R_0^{k}}}|u|\le R_0^{k\alpha}=R_0^{(k+1)\alpha}\frac{1}{R_0^{\alpha}}\le \frac{1}{R_0^{\alpha}}r^{\alpha}.
$$
Therefore for $C(N,p, \mu,\alpha) = \frac{1}{R_0^{\alpha}}$, \eqref{claim step 3.} holds.
\end{proof}
Now we prove that only the smallness in oscillations of coefficients $A_{\pm}$ is sufficient to show the regularity estimates as in Lemma \ref{step 3} for weak solutions of \eqref{PDE}. We prove this result in the following rescaled version of previous lemma.

\begin{lemma}\label{step 3.5}
Suppose $u\in W^{1,p}(B_\rho)$ is bounded and weakly solves the PDE \eqref{PDE} in $B_{\rho}$ and $u(0)=0$.  Then for all $0<\alpha<1$, there exists $C(N,p,\mu,\alpha)>0$ such that for $\delta(N,p, \mu,\alpha)>0$ and $R_0(N,p,\mu, \alpha)$ as in Lemma \ref{step 3} if  
\be\label{SMALLNESS}
\|A_{\pm}-A_{\pm}(0)\|_{L^{\infty}(B_\rho)} < \delta 
\ee
then
\be\label{claim step 3}
\sup_{B_{r}}|u(x)|\le \frac{C(N,p,\mu,\alpha)}{\rho^{\alpha}} \Big (  \Linfty{u}{B_{\rho}} +\rho  \cdot \LN{F}{B_\rho}^{\frac{1}{p-1}} \Big )  r^{\alpha} \;\; \forall r<\rho R_0.
\ee
where $F:=|f_+|+|f_-|$.
\end{lemma}

\begin{proof}
We define the following rescaled function
$$
w(y) := \frac{u(\rho y)}{  \Linfty{u}{B_\rho} + \frac{\rho}{\delta^{\frac{1}{p-1}}}\LN{F}{B_\rho}^{\frac{1}{p-1}}},\;\;y\in B_1.
$$
We can easily verify that 
\be\label{linfty}
\Linfty{w}{B_1}\le 1.
\ee
We can also check from Remark \ref{rescaling 0.1} that $w$ is a weak solution of the following PDE 
\be\label{Jtilde}
\dive(\bar A(x,w)|\nabla w|^{p-2} \nabla  w) = \bar f(x,w) \;\;\mbox{in $B_1$}
\ee
where $\bar A_{\pm}, \bar f_{\pm}$ are defined as follows
\[
\begin{split}
\bar A_{\pm}(y)&:= A_{\pm}(\rho y),\\
\bar f_{\pm}(y)&: = \frac{\rho^p\cdot f_{\pm}(\rho y)}{\left ( \Linfty{u}{B_{\rho }}+\frac{\rho}{\delta^{\frac{1}{p-1}}}\LN{F}{B_{\rho}}^{\frac{1}{p-1}} \right )^{p-1}},\\
\end{split}
\]
We claim that the function $w$ satisfies the assumptions of Lemma \ref{step 3}. Indeed, we can see that from \eqref{SMALLNESS} that 
\be\label{coeff1.1}
\Linfty{\bar A_{\pm} - \bar A_{\pm}(0)}{B_1}\le \delta.
\ee
also since $p>1$ we have 
\be\label{coeff2.1}
\LN{\bar f_{\pm}}{B_1} = \frac{\rho^{p-1}\cdot \LN{f_{\pm}}{B_{\rho}(x_0)}}{\left ( \Linfty{u}{B_{\rho}(x_0)}+\frac{\rho}{\delta^{\frac{1}{p-1}}}\LN{F}{B_{\rho}(x_0)}^{\frac{1}{p-1}} \right )^{p-1}}\le \delta.
\ee
Therefore from Lemma \ref{step 3} we have 
\be\label{estonw}
\sup_{B_{r}}|w(x)|\le C(N,p, \mu,\alpha) r^{\alpha} \;\; \forall r<R_0.
\ee
On rescaling $w$ back to $u$ we obtain 
$$
\sup_{B_{r}}|u(x)|\le \frac{C(N,p, \mu,\alpha)}{\rho^{\alpha}} \left (\Linfty{u}{B_{\rho}} +\rho \LN{f}{B_\rho}^{\frac{1}{p-1}} \right )  r^{\alpha} \;\; \forall r<\rho R_0.
$$
\end{proof}

\begin{remark}\label{rem3.4}
The PDE \eqref{PDE} satisfy the structural condition in \cite[Chapter 10, Section 1]{ED09} and therefore a weak solution $u$ belongs to De-Giorgi class $DG_p(\mu, N) $. Therefore, it is locally Hölder continuous in $B_1$( c.f. \cite[Chapter 10, Theorem 3.1]{ED09}). This means, the sets $\{u>0\}$ and $\{u<0\}$ are open sets. By considering test functions supported inside $\{u>0\}$ or $\{u<0\}$, we observe that any weak solution to \eqref{PDE} is also weak solution to following two PDEs
$$
\begin{cases}
-\dive(A_+(x)|\nabla u|^{p-2}\nabla u)=f_+\qquad \mbox{in $\{u>0\}\cap B_{1}$}\\
-\dive(A_-(x)|\nabla u|^{p-2}\nabla u)=f_-\qquad \mbox{in $\{u<0\}\cap B_1$}.
\end{cases}
$$
From the standard elliptic regularity theory we know, for any given $0<\alpha<1$, $u$ is locally $C^{0,\alpha}$ in $\left ( \{u>0\}\cup \{u<0\} \right )\cap B_{1}$ (c.f. \cite[Theorem 2.4, Section 4]{MS22}). 

We assume $u$ is a bounded weak solution of PDE \eqref{PDE} in $B_1$. From any ball $B_r(x_0)\subset \subset (\{u>0\}\cup \{u<0\})$ we have following estimate (c.f. Appendix \ref{Appendix A}) (here $F:= |f_+|+|f_-|$)
\be\label{PDE estimate}
\|u\|_{C^{\alpha}(B_r(x_0))}\le \frac{1}{r^\alpha} \Bigg (\|u\|_{L^{\infty}(B_1)} +\L{N}{F}{B_1}^{\frac{1}{p-1}} \Bigg )
\ee
However, the above regularity estimates on $u$ deteriorate as we move closer to the free boundary $\partial \{u\neq 0\}\cap B_1$ (since $r\to 0$ as we move close to the free boundary). Therefore, we cannot yet conclude that $u\in C_{loc}^{0,\alpha}(B_{1})$. In order to prove it, we utilize the non-homogenous Moser-Harnack inequality along with some localized geometric arguments.
\end{remark}

\begin{lemma}\label{step 4}
Suppose $u\in W^{1,p}(B_1)$ is a bounded weak solution of \eqref{PDE} in $B_1$.  Then for every $0<\alpha<1$, there exists $\delta(N,p, \mu,\alpha)>0$ such that if 
$$
\|A_{\pm}-A_{\pm}(0)\|_{L^{\infty}(B_1)} < \delta 
$$
then we have 
$$
\|u\|_{C^{\alpha}(B_{1/2})}\le C(N,p, \alpha, \mu) \Big (\Linfty{u}{B_1}+\LN{F}{B_{1}}^{\frac{1}{p-1}} \Big ).
$$
where $F:= |f_+|+|f_-|$.
\end{lemma}

\begin{proof}
In this proof,  $R_0$ is as in Proposition \ref{step 3.5} and $F:=|f_+|+|f_-|$. For the sake of this proof we introduce the following function in the set $\{u\neq 0\}\cap B_1$
$$
d(x) = \begin{cases}
\dist (x, \overline{ \{u\le 0\}})\;\;\mbox{if $u(x)>0$}\\
\dist (x,\overline{ \{u\ge 0\}})\;\;\mbox{if $u(x)<0$}.
\end{cases}
$$
We start by proving the following auxiliary estimates
\begin{enumerate}[label=\textbf{(e-\arabic*)}]
\item \label{e1} For any $y\in {B_{1/2}}$ and $x\in \{u=0\}\cap{ B_{5/8}}$ we have
\be\label{e1eq}
|u(x)-u(y)| = |u(y)| \le C_1(N,p,\mu,\alpha ) \Big ( \Linfty{u}{B_1}+\LN{F}{B_1}^{\frac{1}{p-1}} \Big )|x-y|^{\alpha}.
\ee
\item \label{e1.5}
For any $x\in {B_{1/2}}$ 
\be\label{e1.5eq}
|u(x)| \le C_1(N,p,\mu,\alpha)\Big ( \Linfty{u}{B_1}+\LN{F}{B_1}^{\frac{1}{p-1}} \Big ) d(x)^{\alpha}.
\ee
\item \label{e2} For any $x\in  \Big ( \{u>0\} \cup \{u<0\} \Big ) \cap B_{1/2}$ such that $d:=d(x)\le \frac{R_0}{8}$ 
\be\label{e2eq}
\|u\|_{C^{0,\alpha}(B_{d/8}(x))} \le \frac{C_2(N,p,\mu,\alpha)}{d^{\alpha}} \Big ( u(x) + d \cdot \LN{F}{B_1}^{\frac{1}{p-1}} \Big ).
\ee
\item \label{e3} For any $x\in \Big ( \{u>0\} \cup \{u<0\} \Big ) \cap B_{1/2}$ such that $d=d(x)\le \frac{R_0}{8}$
\be\label{e3eq}
\|u\|_{C^{0,\alpha}(B_{d/8}(x))} \le C_3(N,p,\mu, \alpha) \Big (  \Linfty{u}{B_1}+\LN{F}{B_1}^{\frac{1}{p-1}} \Big ).
\ee
\end{enumerate}
Before delving into the proofs of \ref{e1}-\ref{e3}. We start by observing that for any $x\in B_{1}$ we have
\be\label{smallness around x0}
\Linfty{A_{\pm}-A_{\pm}(x)}{B_{d}(x)} \le \Linfty{A_{\pm}-A_{\pm}(0)}{B_{d}(x)} + |A_{\pm}(0)-A_{\pm}(x)| \le \frac{\delta}{2} +\frac{\delta}{2} =\delta.
\ee
In order to prove \ref{e1}, we observe that $B_{1/4}(x)\subset \subset B_1$. Once $u$ is a weak solution of \eqref{PDE} in $B_1$, so it is in $B_{1/4}$. We now divide the proof of \ref{e1} in two cases
\\[3pt]
\textbf {Case e1.A} : $y\in B_{1/2}$, $x\in B_{5/8}$ and $|x-y|< \frac{R_0}{4}$. 

\vspace{0.3cm}

In this case, since \eqref{smallness around x0} holds, the choice of $\rho=\frac{1}{4}$, $r=|x-y|$ and $x_0=x$ is an admissible choice in Proposition \ref{step 3.5}. It readily follows that for some constant $C_0:= C_0(N,p,\mu,\alpha)$ 
$$
|u(x)-u(y)| = |u(y)| \le C_0(N,p,\mu,\alpha ) \Big ( \Linfty{u}{B_1}+\LN{F}{B_1}^{\frac{1}{p-1}} \Big )|x-y|^{\alpha}.
$$
\textbf{Case e1.B}: $y\in B_{1/2}$, $x\in B_{5/8}$ and $|x-y|\ge \frac{R_0}{4}$. 

$$
{|u(x)-u(y)|} \le \frac{4^\alpha\cdot 2\Linfty{u}{B_1}}{R_0^{\alpha}}|x-y|^{\alpha}.
$$
\noindent
Now, from the \textbf{Case e1.A} and \textbf{Case e1.B}, for $y\in B_{1/2}$ and $x\in \{u=0\}\cap B_{5/8}$ we have
$$
|u(x)-u(y)| \le \max \Big ( C_0, \frac{2\cdot 4^\alpha}{R_0^\alpha} \Big ) \Big ( \Linfty{u}{B_1} + \LN{F}{B_1}^{\frac{1}{p-1}}\Big )|x-y|^{\alpha}.
$$
This proves \ref{e1} with $C_1 = \max \big ( C_0, \frac{2\cdot 4^\alpha}{R_0^\alpha}\big )$. 

In order to prove \ref{e1.5}, let us take $\bar x \in \{u=0\}$ be such that  $d=d(x) = |x-\bar x|$. We again divide the proof in two cases
\\[5pt]
\textbf{Case e2.A}: Assume $d(x)\le \frac{R_0}{8}$.

We observe that $\bar x\in B_{5/8}\cap \{u=0\}$, indeed,
$$
|\bar x| \le |x| + |x-\bar x| \le \frac{1}{2}+\frac{R_0}{8} < \frac{5}{8}.
$$
Also, we can easily verify that $B_{1/4}(\bar x) \subset \subset B_1$.  Hence using \ref{e1} with $\bar  x \in \{u=0\}\cap B_{5/8}$, we have 
\[
\begin{split}
|u(x)-u(\bar x)| = |u(x)| &\le C_{00}(N,p,\mu,\alpha ) \Big ( \Linfty{u}{B_1}+\LN{F}{B_1}^{\frac{1}{p-1}} \Big )|x-\bar x|^{\alpha}\\
&= C_{00}(N,p,\mu,\alpha ) \Big ( \Linfty{u}{B_1}+\LN{F}{B_1}^{\frac{1}{p-1}} \Big )d(x)^{\alpha}.
\end{split}
\]
\\
\textbf{Case e2.B}: Assume that $d(x)>\frac{R_0}{8}$. In this case 
$$
|u(x)| \le \frac{8^\alpha \cdot \Linfty{u}{B_1}}{R_0^\alpha}d(x)^{\alpha}.
$$
\noindent
Thus, \textbf{Case e2.A} and \textbf{Case e2.B} prove \ref{e1.5} with $C_2 = \max \Big ( C_{00},  \frac{4^\alpha }{R_0^\alpha}\Big )$. 
\\[5pt]
For the proof of \ref{e2}, it is enough to consider only the case where $x\in B_{1/2}\cap \{u>0\}$, since the case where $x\in B_{1/2}\cap \{u<0\}$ can be treated similarly. We know that $u$ is a weak solution of the following PDE
\be\label{PDE++}
\dive\Big ( A_+(x)|\nabla u|^{p-2}\nabla u \Big ) = f_+ \;\;\mbox{in $B_{d/4}(x)$}.
\ee
From the  non-homogenous Moser-Harnack inequality \cite[Theorem 1]{serrin63}, we have  
\be\label{harnack.}
\begin{split}
\sup_{B_{d/4}(x)} u& \le C (N,p,\mu) \Big  (  \inf_{B_{d/8}(x)}u+d \cdot \|f_+\|_{L^N(B_{d/4}(x))}^{\frac{1}{p-1}} \Big  )\\
&\le C (N,p,\mu) \Big  (  \inf_{B_{d/8}(x)}u+d \cdot \|F\|_{L^N(B_{d/4}(x))}^{\frac{1}{p-1}} \Big  ).
\end{split}
\ee
Moreover from the regularity estimates for $u$ in $\{u>0\}$ we have
\be\label{schauder and harnack}
\begin{split}
\|u\|_{C^{0,\alpha}(B_{d/8}(x))} &\le \frac{C(N,p,\mu,\alpha)}{d^{\alpha}}  \left ( \sup_{B_{d/4}}(u)+d\cdot \LN{F}{B_{d/4}(x)}^{\frac{1}{p-1}} \right ).
\end{split}
\ee
Using \eqref{harnack.} and \eqref{schauder and harnack} we arrive at
\[
\begin{split}
\|u\|_{C^{0,\alpha}(B_{d/8}(x))} &\le \frac{C (N,p,\mu,\alpha)}{d^{\alpha}} \left (  \inf_{B_{d/8}(x)}u+d \cdot \|F\|_{L^N(B_{d/4}(x))} ^{\frac{1}{p-1}}\right )\\
&\le \frac{C (N,p,\mu,\alpha)}{d^{\alpha}} \left (  u(x)+d \cdot \|F\|_{L^N(B_{d/4}(x))}^{\frac{1}{p-1}} \right ).
\end{split}
\]
This concludes the proof of \ref{e2}. In order to prove \ref{e3}, we again treat only the case $x\in \{u>0\}\cap B_{1/2}$. From \ref{e1.5}
$$
u(x)=|u(x)| \le C(N,p,\mu ,\alpha)\Big ( \Linfty{u}{B_1}+\LN{F}{B_1} ^{\frac{1}{p-1}}\Big ) d(x)^{\alpha}.
$$
Plugging the above estimates in \ref{e2} we obtain 
\be\label{+phase}
\begin{split}
\|u\|_{C^{\alpha}(B_{d(x)/8}(x))}&\le C_1C_2 \left (\Linfty{u}{B_1}+\LN{F}{B_{1}}^{\frac{1}{p-1}} \right )+C_1d(x)^{1-\alpha}\|F\|_{L^N(B_1)}^{\frac{1}{p-1}}\\
&\le C_3(N, p,\mu, \alpha) \left (\Linfty{u}{B_1}+\LN{F}{B_{1}}^{\frac{1}{p-1}} \right ).
\end{split}
\ee
Now, under the possession of \ref{e1}-\ref{e3}, we now finish the proof of Lemma \ref{step 4}. 

We again divide the proof in cases.
\\[5pt]
\textbf{Case I}.  Let $x,y\in B_{1/2}$ are such that $u(x)\cdot u(y)=0$. 
\\[5pt]
We can assume without loosing generality that $u(x)=0$. Then, it follows readily from \ref{e1} that 
$$
|u(x)-u(y)| = |u(y)| \le C_1(N,p,\mu,\alpha ) \Big ( \Linfty{u}{B_1}+\LN{F}{B_1}^{\frac{1}{p-1}} \Big )|x-y|^{\alpha}.
$$
\textbf{Case II}. Let $x,y\in B_{1/2}$ such that $u(x)\cdot u(y)\neq 0$. 
\\[5pt]
Without loss of generality, we can assume the 
$$
d(y)\le d(x).
$$
Once more, splitting the proof in cases
 
\vspace{0.2cm}
\noindent
\textbf{Case II.1}. If $|x-y|< \frac{d(x)}{8}$.  We now study the two subcases

\vspace{0.3cm}
\hspace{-0.5cm}
\textbf{Case II.1.A}. If $d(x)\le \frac{R_0}{8}$. 

\indent\indent
\hspace{-0.7cm}{In this case, $y\in B_{d(x)/8}(x)$.} 			Then it readily follows from \ref{e3} that 
			\[
			\begin{split}
			|u(x)-u(y)| &\le [u]_{C^{\alpha}(B_{d(x)/8}(x))}|x-y|^{\alpha}
			\le C_3 \Big ( \Linfty{u}{B_1}+\LN{F}{B_1}^{\frac{1}{p-1}} \Big )|x-y|^{\alpha}
			\end{split}
			\]
\indent \indent where $C_3:=C_3(N,p,\mu,\alpha)$.	
\vspace{0.3cm}
		
\hspace{-0.5cm} \textbf{Case II.1.B}. If $d(x) > \frac{R_0}{8}$. 

\indent\indent
Since $\dive\big ( A_+(x)|\nabla u|^{p-2}\nabla u \big ) = f_+ \;\;\mbox{in $B_{d/8}(x)$}$ in weak sense. This leads to 
		
			\[
			\begin{split}
			|u(x)-u(y)| &\le [u]_{C^{\alpha}(B_{d(x)/8}(x))}|x-y|^{\alpha}\\
			&\le \frac{C(N,p,\mu,\alpha)}{d^{\alpha}} \Big ( \Linfty{u}{B_1}+d \LN{F}{B_1}^{\frac{1}{p-1}} \Big )|x-y|^{\alpha}\\
			&\le \frac{C_4(N,p,\mu,\alpha)}{R_0^{\alpha}} \Big (  \Linfty{u}{B_1}+\LN{F}{B_1}^{\frac{1}{p-1}} \Big )|x-y|^{\alpha}.
			\end{split}
			\]

\textbf{Case II.2}. If $|x-y|\ge  \frac{d(x)}{8}$. 
\\[5pt]
\indent \indent By \ref{e1.5} and the assumption $d(x)\ge d(y)$ we obtain 
\be\label{more}
\begin{split}
|u(x)-u(y)| &= |u(x)|+|u(y)|\\
&\le  C_1(N, p,\mu, \alpha) \Big (\Linfty{u}{B_1}+\LN{F}{B_{1}}^{\frac{1}{p-1}}\Big ) (d(x)^{\alpha}+d(y)^{\alpha})\\
&\le 2  C_1(N,p, \mu,\alpha) \Big (\Linfty{u}{B_1}+\LN{F}{B_{1}}^{\frac{1}{p-1}}\Big )d(x)^{\alpha}\\
&\le C_5(N, p,\mu,\alpha) \Big (\Linfty{u}{B_1}+\LN{F}{B_{1}}^{\frac{1}{p-1}}\Big )|x-y|^{\alpha}.
\end{split}
\ee
\noindent This proves Lemma \ref{step 4}.
\end{proof}

We can rescale the above lemma to a ball of any radius $\rho$.
\begin{corollary}\label{step 4.5}
Let $u\in W^{1,p}(B_\rho(x_0))$ be a bounded weak solution of \eqref{PDE} in $B_{\rho}(x_0)$ and $\rho< 1$.  Then for every $0<\alpha<1$ and for $\delta(N, p,\mu,\alpha)>0$ as in Lemma \ref{step 3}, we have 
\be\label{claim 4}
\|A_{\pm}-A_{\pm}(x_0)\|_{L^{\infty}(B_\rho(x_0))} \le \frac{\delta}{2}  \implies \|u\|_{C^{\alpha}(B_{\frac{\rho}{2}}(x_0))} \le \frac{C}{\rho^{\alpha}}\left (\Linfty{u}{B_\rho(x_0)}+\rho \LN{F}{B_\rho(x_0)}^{\frac{1}{p-1}}\right )
\ee
where $C:=C(N, p,\mu, \alpha)$.
\end{corollary}
\begin{proof}
We reduce the Corollary \ref{step 4.5} to Proposition \ref{step 4} by using  the following rescaling 
$$
w(y) = \frac{1}{\rho}u(x_0+\rho y).
$$
From Remark \ref{rescaling 0.1}, we prove the  Corollary \ref{step 4.5}.
\end{proof}

\section{Proof of Theorem \ref{main result}}

And now, present the proof of the Theorem \ref{main result}.

\begin{proof}[Proof of Theorem \ref{main result}]

Let $F:=|f_+|+|f_-|$ as before and $D\subset \subset B_1$. We observe that $D\subset \subset B_{d} \subset \subset B_1$ where $d=1-\frac{\dist(D,\partial B_1)}{2}$. Since $A_{\pm}\in C(B_1)$, then $A_{\pm}$ are uniformly continuous in $\overline{B_d}$. Thus, we define $\o_{A_{\pm}, B_d}$ as the following  modulus of continuity. We define $\o_{A_{\pm},\overline{B_d}}$ to be 
$$
\o_{A_{\pm},\overline{B_d}}(t) := \max \Bigg ( \sup_{\substack{ |x-y|<t\\x,y\in \overline{B_d}}}|A_+(x)-A_+(y)|\,,\,\sup_{\substack{ |x-y|<t\\x,y\in \overline{B_d}}}|A_-(x)-A_-(y)|\Bigg )\mbox{  for $t\le 2d$}
$$
and we define $\o_{A_{\pm},\overline{B_d}}(t) := \o_{A_{\pm},\overline{B_d}}(2d)$ for $t>2d$.
{We set $t_0$ as  
$$
t_0 := t_0(\o_{A_{\pm},\overline{B_d}},\delta) = \sup \Big \{ t \;\Big | \;\o_{A_{\pm},\overline{B_d}}(t)\le \delta \Big \}
$$ 
as well as 
$$
s_0 := \min\Bigg (t_0,\frac{\dist(D,\partial B_1)}{4} \Bigg ).
$$
Since $\o_{A_{\pm},\overline {B_d}}$ is a non-decreasing function we have $\o_{A_{\pm},\overline{B_d}}(s_0)\le \delta$. }Now since 
$$
D\subset \bigcup_{x\in D}B_{s_0}(x) \subset \overline {B_d}
$$
we have 
\be\label{smallness1...}
\sup_{B_{s_0}(x)}|A_{\pm}-A_{\pm}(x_0)| \le \o_{A_{\pm},\overline{B_d}}(s_0) \le {\delta},\;\forall x\in D.
\ee
$u$ is a weak solution of \eqref{PDE} in $B_{s_0}(x)$, we have from Corollary \ref{step 4.5} that for all $y\in B_{s_0/2}(x)\cap D$ 
\be\label{near x_0...}
|u(x)-u(y)|\le \frac{C(N,p, \mu,\alpha)}{s_0^{\alpha}}\left (\Linfty{u}{B_{1}}+s_0\LN{F}{B_{1}}^{\frac{1}{p-1}}  \right )|x-y|^{\alpha}.
\ee
Now, if $x, y\in D$ are such that $|x-y| \ge s_0/2$, then 
\be\label{far from x_0...}
{|u(x)-u(y)|}\le 2^{1+\alpha} \frac{\Linfty{u}{B_1}}{s_0^{\alpha}}|x-y|^{\alpha}.
\ee
By combining \eqref{near x_0...} and \eqref{far from x_0...}, we arrive to $\Big ($since $s_0\le\frac{\dist(D,\partial B_1)}{4}\le\frac{\diam(B_1)}{4} =\frac{1}{2}<1$$\Big )$
\be\label{precise est...}
[u]_{C^{\alpha}(D)} \le \frac{C(N,p,\alpha, \mu)}{s_0^{\alpha}}\left ( \Linfty{u}{B_{1}}+\LN{F}{B_{1}} ^{\frac{1}{p-1}}\right ).
\ee
{We observe from the definition of $s_0$ that 
$$
[u]_{C^{\alpha}(D)} \le \begin{cases}
\frac{C(N, p,\mu,\alpha)}{t_0^{\alpha}}\left ( \Linfty{u}{B_{1}}+\LN{F}{B_{1}}^{\frac{1}{p-1}} \right ) \;\;\;\;\;\;\;\;\;\;\;\;\;\;\;\;\;\;\;\mbox{if $t_0\le \frac{\dist(D,\partial B_1)}{4}$}\\[7pt]
\frac{4^\alpha\cdot C(N, p,\mu,\alpha)}{\dist(D,\partial B_1)^\alpha}\left ( \Linfty{u}{B_{1}}+\LN{F}{B_{1}}^{\frac{1}{p-1}} \right )\;\; \;\;\;\;\;\;\;\;\;\;\;\;\;\mbox{if $t_0\ge \frac{\dist(D,\partial B_1)}{4}$.}
\end{cases}
$$
}
{In order to control the first term in the equation above by a universal multiple of $\dist(D,\partial B_1)^{-\alpha}$, we observe that} once $t_0>0$ depends only on the modulus of continuity $\o_{A_{\pm}, \overline{B_d}}$ and $\delta$, there exists a $n_0:= n_0(\o_{A_{\pm}, \overline{B_d}},\delta)=n_0(N,p,\mu,\alpha,\o_{A_{\pm},\overline{B_d}})\in \N$ such that  $\frac{2}{n_0}\le t_0$. Hence
$$
\frac{\dist(D,\partial B_1)}{n_0}\le \frac{\diam(B_1)}{n_0} =\frac{2}{n_0}\le t_0 \implies \frac{1}{t_0^{\alpha}}\le \frac{n_0^{\alpha} }{\dist(D,\partial B_1)^{\alpha}}= \frac{C(N,p,\mu,\alpha,\o_{A_{\pm},\overline{B_d}}) }{\dist(D,\partial B_1)^{\alpha}}.
$$
Now, \eqref{precise est...} becomes
$$ 
[u]_{C^{\alpha}(D)} \le \frac{C(N,p,\mu,\alpha,\o_{A_{\pm},\overline{B_d}})}{\dist(D,\partial B_1)^{\alpha}}\left ( \Linfty{u}{B_{1}}+\LN{F}{B_{1}}^{\frac{1}{p-1}} \right ).
$$
By observing that 
$$
\L{N}{F}{B_1}^{\frac{1}{p-1}} \le  \Big ( \L{N}{f_+}{B_1} + \L{N}{f_-}{B_1} \Big )^{\frac{1}{p-1}} \le  C(p)\Big ( \L{N}{f_+}{B_1}^{\frac{1}{p-1}} + \L{N}{f_-}{B_1}^{\frac{1}{p-1}} \Big )
$$
we prove the estimate \eqref{general est...}. Finally, in order to obtain \eqref{estimates}, we just take $D:= B_r$.
\end{proof}

\begin{appendices}
\section{Proof of Theorem \ref{supporting lemma}}\label{Appendix A}
Since the PDE \eqref{simple PDE} is invariant under translation and addition of a constant, we can assume $x_0=0\in \O$ and $u(0)=0$ without loosing generality. The proof of the following lemma is exactly the same as the proof of Lemma \ref{step 3} , with $A_{\pm}=A$ and $f_{\pm}=f$. 
\begin{lemma}\label{step .3}
Suppose $u\in W^{1,p}(B_1)$ weakly satisfy the PDE \eqref{simple PDE} with $\| u\|_{L^{\infty}(B_1)}\le 1$  and $u(0)=0$.  Then for all $0<\alpha<1$ and $\delta(N,p, \mu,\alpha)>0$ as in Lemma \ref{step 2} $($with $A_{\pm}=A$ and $f_{\pm}=f$$)$  there exists $C(N,p,\mu,\alpha)>0$ such that if 
$$
 \max \Big (\|A-A(0)\|_{L^{\infty}(B_1)},\, \|f\|_{L^N(B_1)} \Big ) < \delta 
$$
then for $R_0(N,p,\mu, \alpha)$ as in Lemma \ref{step 2} $($with $A_{\pm}=A$ and $f_{\pm}=f$$)$ we have
\be\label{claim step 3..}
\sup_{B_{r}}|u(x)|\le C(N,p, \mu,\alpha) \cdot r^{\alpha} \;\; \forall r<R_0.
\ee
Precisely speaking, we have $C(N,p,\mu,\alpha)= R_0^{-\alpha}$.
\end{lemma}
Just like we proved the Lemma \ref{step 3.5}, we prove the following Lemma \ref{step .3.5.} by taking $A_{\pm}=A$ and $f_{\pm}=f$.

\begin{lemma}\label{step .3.5.}
Suppose $u\in W^{1,p}(B_\rho)$ be bounded and weakly solves the PDE \eqref{simple PDE} in $B_{\rho}$ and $u(0)=0$.  Then for all $0<\alpha<1$, there exists $C(N,p,\mu,\alpha)>0$ such that for $\delta(N,p, \mu,\alpha)>0$ and $R_0(N,p, \mu,\alpha)$ as in Lemma \ref{step .3} if  
\be\label{SMALLNESS.}
\|A-A(0)\|_{L^{\infty}(B_\rho)} < \delta 
\ee
then
\be\label{claim step 3..}
\sup_{B_{r}}|u(x)|\le \frac{C(N,p,\mu,\alpha)}{\rho^{\alpha}} \Big  (  \Linfty{u}{B_{\rho}} +\rho \cdot \LN{f}{B_\rho}^{\frac{1}{p-1}} \Big ) \cdot  r^{\alpha} \;\; \forall r<\rho R_0.
\ee
\end{lemma}

The proof of Theorem \ref{supporting lemma} readily follows from the rescaling and covering argument as in the proof of Theorem \ref{main result} by assuming $A_+ = A_-$ and $f_+ =f_-$. 

\end{appendices}
\section*{Acknowledgements}
Authors thank Prof. Peter Lindqvist and Prof. Juan Manfredi for their insightful remarks during preparation the paper. Authors also thank Institut Mittag-Leffler, Djursholm (Sweden) for the invitation to participate on the event \textbf{``Geometric aspects of nonlinear partial differential equations" } in the fall of 2022, where the discussions about this paper took place.  The first author thanks the financial support  from CNPq-311566/2019-7 (Brazil) and FUNCAP (PRONEX). The second author is also thankful to the financial support from TIFR-CAM, Bangalore (India).

 \newcommand{\noop}[1]{}

\end{document}